\documentclass[a4paper,11pt]{article}

\makeatletter
 
  \@addtoreset{equation}{section}
 \makeatother
\usepackage{amsfonts}
\usepackage{amssymb}
\usepackage{mathrsfs}
\usepackage{amsmath,amsthm}
\usepackage{amsmath}
\usepackage{comment}
\usepackage{algorithm}
\usepackage{algorithmic}

\usepackage[dvipdfmx]{graphicx}
\usepackage{color}
\usepackage{indentfirst}
\usepackage{here}
\usepackage{url}
\usepackage[hang,small,bf]{caption}
\usepackage[labelformat=parens]{subcaption}
\begin{document}

\def\bn{{\bf n}}
\def\A{{\bf A}}
\def\B{{\bf B}}
\def\C{{\bf C}}
\def\D{{\bf D}}
\def\E{{\bf E}}
\def\F{{\bf F}}
\def\G{{\bf G}}
\def\H{{\bf H}}
\def\I{{\bf I}}
\def\J{{\bf J}}
\def\K{{\bf K}}
\def\L{{\bf L}}
\def\M{{\bf M}}
\def\N{{\bf N}}
\def\O{{\bf O}}
\def\P{{\bf P}}
\def\Q{{\bf Q}}
\def\R{{\bf R}}
\def\S{{\bf S}}
\def\T{{\bf T}}
\def\U{{\bf U}}
\def\V{{\bf V}}
\def\W{{\bf W}}
\def\X{{\bf X}}
\def\Y{{\bf Y}}
\def\Z{{\bf Z}}
\def\cala{{\cal A}}
\def\calb{{\cal B}}
\def\calc{{\cal C}}
\def\cald{{\cal D}}
\def\cale{{\cal E}}
\def\calf{{\cal F}}
\def\calg{{\cal G}}
\def\calh{{\cal H}}
\def\cali{{\cal I}}
\def\calj{{\cal J}}
\def\calk{{\cal K}}
\def\call{{\cal L}}
\def\calm{{\cal M}}
\def\caln{{\cal N}}
\def\calo{{\cal O}}
\def\calp{{\cal P}}
\def\calq{{\cal Q}}
\def\calr{{\cal R}}
\def\cals{{\cal S}}
\def\calt{{\cal T}}
\def\calu{{\cal U}}
\def\calv{{\cal V}}
\def\calw{{\cal W}}
\def\calx{{\cal X}}
\def\caly{{\cal Y}}
\def\calz{{\cal Z}}
%
\def\sskip{\hspace{0.5cm}}
\def\simleq{ \raisebox{-.7ex}{\em $\stackrel{{\textstyle <}}{\sim}$} }
\def\leqsim{ \raisebox{-.7ex}{\em $\stackrel{{\textstyle <}}{\sim}$} }
\def\ep{\epsilon}
\def\half{\frac{1}{2}}
\def\iku{\rightarrow}
\def\Iku{\Rightarrow}
\def\ikup{\rightarrow^{p}}
\def\inclusion{\hookrightarrow}
\def\cadlag{c\`adl\`ag\ }
\def\up{\uparrow}
\def\down{\downarrow}
\def\doti{\Leftrightarrow}
\def\douti{\Leftrightarrow}
\def\dochi{\Leftrightarrow}
\def\douchi{\Leftrightarrow}%
\def\yy{\\ && \nonumber \\}
\def\y{\vspace*{3mm}\\}
\def\nn{\nonumber}
\def\be{\begin{equation}}
\def\ee{\end{equation}}
\def\bea{\begin{eqnarray}}
\def\eea{\end{eqnarray}}
\def\beas{\begin{eqnarray*}}
\def\eeas{\end{eqnarray*}}
%
\def\hd{\hat{D}}
\def\hv{\hat{V}}
\def\hsd{{\hat{d}}}
\def\hx{\hat{X}}
\def\hsx{\hat{x}}
\def\bsx{\bar{x}}
\def\bsd{{\bar{d}}}
\def\bx{\bar{X}}
\def\ba{\bar{A}}
\def\bb{\bar{B}}
\def\bc{\bar{C}}
\def\bv{\bar{V}}
\def\balpha{\bar{\alpha}}
\def\bbalpha{\bar{\bar{\alpha}}}
\def\combi{\l(\begin{array}{c}\alpha\\ \beta \end{array}\r)}
\def\f{^{(1)}}
\def\s{^{(2)}}
\def\ss{^{(2)*}}
\def\l{\left}
\def\r{\right}
\def\a{\alpha}
\def\b{\beta}
\def\L{\Lambda}

\newtheorem{thm}{Theorem}
\newtheorem{lemma}{Lemma}
\newtheorem{prop}{Proposition}
\newtheorem{defn}{Definition}
\newtheorem{rem}{Remark}
\newtheorem{step}{Step}
\newtheorem{cor}{Corollary}
\newtheorem{assumption}{Assumption}
\newtheorem{ex}{Example}

\everymath {\displaystyle}

\newcommand{\ruby}[2]{
\leavevmode
\setbox0=\hbox{#1}
\setbox1=\hbox{\tiny #2}
\ifdim\wd0>\wd1 \dimen0=\wd0 \else \dimen0=\wd1 \fi
\hbox{
\kanjiskip=0pt plus 2fil
\xkanjiskip=0pt plus 2fil
\vbox{
\hbox to \dimen0{
\small \hfil#2\hfil}
\nointerlineskip
\hbox to \dimen0{\mathstrut\hfil#1\hfil}}}}

\def\qedsymbol{$\blacksquare$}
\renewcommand{\thefootnote }{\fnsymbol{footnote}}
\renewcommand{\refname }{References}

\everymath {\displaystyle}

\title{
New asymptotic expansion formula via Malliavin calculus and its application to rough differential equation driven by fractional Brownian motion
}
\author{Akihiko Takahashi\footnote{University of Tokyo, Japan} \ and  Toshihiro Yamada\footnote{Hitotsubashi University, Japan} \footnote{Corresponding author. Phone: +81-42-580-8788, 
E-mail: toshihiro.yamada@r.hit-u.ac.jp}} 
\date{ }
\maketitle

\abstract
This paper presents a novel generic asymptotic expansion formula of expectations of multidimensional Wiener functionals through a Malliavin calculus technique. The uniform estimate of the asymptotic expansion is shown under a weaker condition on the Malliavin covariance matrix of the target Wiener functional. In particular, the method provides a tractable expansion for the expectation of an irregular functional of the solution to a multidimensional rough differential equation driven by fractional Brownian motion with Hurst index $H<1/2$, without using complicated fractional integral calculus for the singular kernel. In a numerical experiment, our expansion shows a much better approximation for a probability distribution function than its normal approximation, which demonstrates the validity of the proposed method. \\
\\
{\bf Keywords}:  Asymptotic expansion, Wiener functional, Malliavin calculus, Rough differential equation, Fractional Brownian motion

\section{Introduction}

In the paper, we derive a new asymptotic expansion formula of the expectations of multidimensional Wiener functionals as an extension of Watanabe (1987) \cite{W87}, Yoshida (1992) \cite{Yo}, Takahashi (1999) \cite{T1}, Kunitomo and Takahashi (2001, 2003) \cite{KT1,KT2}, Takahashi and Yoshida (2005) \cite{TYo}, Malliavin and Thalmaier (2006) \cite{MallThal}, Takahashi and Yamada (2012) \cite{TY12} and Takahashi (2015) \cite{T2}.   
 The general asymptotic expansion through a Malliavin calculus approach provides wide applications and covers previous expansion schemes. More precisely, a technique with a Malliavin derivative (annihilation) computation and a Skorohod integral (divergence, creation) computation is introduced. A fractional order expansion on an abstract Wiener space is considered to apply the method to general Gaussian processes, particularly, rough differential equations driven by fractional Brownian motion. The asymptotic expansion of $\mathbb{E}[f(F^\varepsilon)]$ for a multidimensional Wiener functional $F^\varepsilon$ with a small parameter is obtained under a weaker condition in the sense that we only impose an assumption for the inverse of the Malliavin covariance for $F^0$, the dominant part of the expansion $F^\varepsilon \sim F^0+\cdots$, not for $F^\varepsilon$ itself. The condition is always easily checked in practical stochastic models. The test function $f$ is assumed to be a bounded measurable function, and the uniform bound of the expansion is shown.

The method provides a tractable expansion for the expectation of an irregular functional of the solution to a multidimensional rough differential equation driven by fractional Brownian motion with Hurst index $H<1/2$ without using complicated fractional integral calculus for the singular kernel. 
We take an approach substantially different from Baudoin (2015) \cite{Ba}, Inahama (2016) \cite{I} for the asymptotics for the density of a solution to a rough differential equation.  
To obtain explicit expansion formulas for expectations of irregular functionals of solutions to rough differential equations, we use the Stratonovich--Skorohod transformation, the Stroock-Taylor formula and the integration by parts in Malliavin calculus. 
Then the expansion terms involving iterated rough integrals are all transformed into polynomials of fractional Brownian motion which can be easily simulated by Monte Carlo or quasi Monte Carlo methods.  A numerical example for the asymptotic expansion of a probability distribution function is shown to validate the method. In particular, a comparison result with the normal approximation shows the effectiveness of our expansion scheme. 

The paper is organized as follows. Section 2 provides a new asymptotic expansion for general Wiener functionals on an abstract Wiener space. Then, Section 3 shows a tractable expansion for the expectation of an irregular functional of the solution to a multidimensional rough differential equation driven by fractional Brownian motion of the Hurst index $H<1/2$ with a numerical example. Section 4 concludes. 

\section{Asymptotic expansion of expectation of Wiener functionals}

We prepare notation and definitions on Malliavin calculus on an abstract Wiener space. For the details, see 
Ikeda and Watanabe (1989) \cite{IW}, Malliavin (1997) \cite{Mall} and Nualart (2006) \cite{N}. 

Let $({\cal W}, {\cal H}, \mu)$ be an abstract Wiener space, where ${\cal W}$ is a Banach space, ${\cal H}$ is a separable Hilbert space which is continuously, densely embedded into ${\cal W}$ called the Cameron-Martin space. 

For $p\geq 1$ and a Hilbert space $G$ equipped with a norm $\| \cdot \|_G$, let $L^p({\cal W};G)$ be a Banach space of all $\mu$-measurable functionals $F:{\cal W} \to G$ such that $\textstyle{\| X \|_{L^p({\cal W};G)}=(\int_{{\cal W}} \| X(\omega) \|^p_G d\mu(\omega))^{1/p}}<\infty$ with the identification $X = Y$ if and only if $X(\omega)=Y(\omega)$ a.e. If $G=\mathbb{R}$, we may use a notation $L^p({\cal W})=L^p({\cal W};G)$.

We denote by $j:{\cal H} \to {\cal W}$ the embedding map. 
Let ${\cal W}'$ and ${\cal H}'$ be the topological dual spaces of ${\cal W}$ and ${\cal H}$, respectively. Then, ${\cal W}' \overset{j^\ast}{\hookrightarrow} {\cal H}'={\cal H} \overset{j}{\hookrightarrow} {\cal W}$ where $j^\ast$ is the dual map of $j$. The Gaussian measure $\mu$ on $({\cal W},\mathscr{B}(\cal W))$ satisfies 
\begin{align}
\int_{{\cal W}} e^{\mathrm{i} _{{\cal W}'} \langle l,\omega \rangle_{{\cal W}} }  d\mu(\omega)
=
e^{-\frac{1}{2} \| j^\ast(l) \|_{{\cal H}}^2 }, \ \ \mbox{for  all} \ l\in {\cal W}'. 
\end{align} 
Hence $\{ _{{\cal W}'} \langle l,\cdot \rangle_{{\cal W}}; \ l \in {\cal W}' \}$ is a family of Gaussian random variables on $({\cal W},\mathscr{B}(\cal W),\mu)$ with mean $0$ and covariance $\mathbb{E}[_{{\cal W}'} \langle l,\omega \rangle_{{\cal W}} \  _{{\cal W}'} \langle l',\omega \rangle_{{\cal W}}]=\langle j^\ast(l),j^\ast(l') \rangle_{{\cal H}}$, $l,l' \in {\cal W}'$. Thus, the map $j^\ast({\cal W}') \ni j^\ast(l) \ \mapsto \ _{{\cal W}'} \langle l,\cdot  \rangle_{{\cal W}} \in L^2({\cal W})$ is a linear isometry which can be extend to an isometry $I: {\cal H} \to L^2({\cal W})$ such that $I(h)$ is a Gaussian random variable with mean $0$ with the standard deviation $\| I(h) \|_{L^2({\cal W})}=\| h \|_{{\cal H}}$ since $j^\ast({\cal W}) \subset {\cal H}$ is dense.  

Let $\mathscr{S}({\cal W})=\{ F:{\cal W} \to \mathbb{R}; \ F=f(I(h_1),\ldots,I(h_n)), \ n\in \mathbb{N}, \ f \in C_b^\infty(\mathbb{R}^n), h_1,\ldots,h_n \in {\cal H} \}$. For $F=f(I(h_1),\ldots,I(h_n)) \in \mathscr{S}({\cal W})$, we define the Malliavin derivative $DF \in {\cal H}$ as 
\begin{align}
DF=\sum_{i=1}^n (\partial_i f)(I(h_1),\ldots,I(h_n)) h_i.
\end{align}
The operator $D$ is a closable operator, and for $p> 1$, we define $\mathbb{D}^{1,p}=\overline{\mathscr{S}({\cal W})}^{\| \cdot \|_{1,p}}$ where the norm $\| \cdot \|_{1,p}$ given by $\| F \|_{1,p}=\| F \|_{{\cal L}^p({\cal W})}+\| DF \|_{{\cal L}^p({\cal W};{\cal H})}$. Similarly, the higher-order Malliavin derivatives $D^k$ and the corresponding Sobolev spaces $\mathbb{D}^{k,p}$ can be defined iteratively. We define $\mathbb{D}^\infty=\cap_{k\in \mathbb{N},p> 1}\mathbb{D}^{k,p}$. 

For $p \in \mathbb{N}$, let $\mathrm{Dom}\delta^p=\{ u \in L^2({\cal W};{\cal H}^{\otimes p}); \ \exists C>0 \mbox{ s.t. } |\mathbb{E}[\langle D^p F,u\rangle_{{\cal H}^{\otimes p}}]| \leq C \| F \|_{L^2({\cal W})}, \forall F \in \mathbb{D}^{p,2} \}$. For $u \in \mathrm{Dom}\delta^p$, there exists $\delta^p(u) \in L^2({\cal W})$ such that
\begin{align}
\mathbb{E}[\langle D^pF,u\rangle_{{\cal H}^{\otimes p}}]=\mathbb{E}[F \delta^p(u) ],  \label{duality_formula}
\end{align}
which is called the duality formula. 
For $F \in \mathbb{D}^{\infty,2}=\cap_{k\in \mathbb{N}}\mathbb{D}^{k,2}$, we have the {\it Stroock-Taylor formula}: 
\begin{align}
F=\mathbb{E}[F]+\sum_{p\geq 1} \frac{1}{p!} \delta^p(\mathbb{E}[D^pF]). \label{Stroock_formula}
\end{align}
For the Stroock-Taylor formula, see Theorem (6) of Stroock (1987) \cite{Stroock}, Proposition 2 in Chapter IV of  \"{U}st\"{u}nel (1995) \cite{U} and Theorem 4.1 in Section 4, Chapter VI.4 of Malliavin (1997) \cite{Mall}, for instance.

For $F=(F^1,\ldots,F^e) \in (\mathbb{D}^\infty)^e$, we define the Malliavin covariance matrix $\sigma^F=[\sigma^F_{ij}]_{1\leq i,j \leq e}$:  
\begin{align}
\sigma^F_{ij}= \langle DF^i,DF^j \rangle_{{\cal H}}, \ \ \  1\leq i,j \leq e.
\end{align} 
We say $F \in (\mathbb{D}^\infty)^e$ is nondegenerate if $\sigma^F$ is invertible a.s. and 
\begin{align}
\| \det (\sigma^F)^{-1} \|_{{\cal L}^p({\cal W})} < \infty, \ \ \ \forall p> 1.
\end{align} 
Let $\mathcal{S}(\mathbb{R}^e)$ be the Schwartz space or the space of $\mathbb{R}$-valued rapidly decreasing functions on $\mathbb{R}^e$. For a nondegenerate Wiener functional $F \in (\mathbb{D}^\infty)^e$, $G \in \mathbb{D}^\infty$, $f \in \mathcal{S}(\mathbb{R}^e)$ and a multi-index $\alpha \in \{1,\ldots,e \}^k$, we have the integration by parts (IBP) formula: 
\begin{align}
\mathbb{E}[\partial^\alpha f (F)G]=\mathbb{E}[f (F)H_{\alpha}(F,G)], \label{IBP_formula}
\end{align} 
where $H_{\alpha}(F,G)$ is recursively defined by $H_{\alpha}(F,G)=H_{(\alpha_k)}(F,H_{(\alpha_1,\ldots,\alpha_{k-1})}(F,G))$ with 
\begin{align}
&H_{(i)}(F,G)=\sum_{j=1}^e \delta( \gamma^F_{ij} DF^j G), \ \ \ i=1,\ldots,e,
\end{align} 
with the inverse matrix $\gamma^F$ of the Malliavin covariance of $F$, i.e. $\gamma^F=(\sigma^F)^{-1}$. 

Let $\mathcal{S}'(\mathbb{R}^e)$ be the dual of $\mathcal{S}(\mathbb{R}^e)$, i.e. $\mathcal{S}'(\mathbb{R}^e)$ is the space of Schwartz tempered distributions. 
Let $\mathbb{D}^{-\infty}$ be the dual space of $\mathbb{D}^{\infty}$, i.e. the space of continuous linear forms on $\mathbb{D}^{\infty}$. For $T \in {\cal S}'(\mathbb{R}^e)$, a multi-index $\alpha=(\alpha_1,\ldots,\alpha_k)$, a nondegenerate $F \in (\mathbb{D}^\infty)^e$ and $G \in \mathbb{D}^\infty$, we have
\begin{eqnarray}
{}_{\mathbb{D}^{-\infty}} \langle \partial^{\alpha}T(F),G \rangle{}_{\mathbb{D}^\infty}={}_{\mathbb{D}^{-\infty}}\langle T(F),H_{\alpha}(F,G) \rangle {}_{\mathbb{D}^\infty},
\label{ibp_formula}
\end{eqnarray}
where  ${}_{{\cal S}'}\langle \cdot, \cdot \rangle_{{\cal S}} $ is the bilinear form on ${\cal S}'(\mathbb{R}^e)$ and ${\cal S}(\mathbb{R}^e)$,  ${}_{\mathbb{D}^{-\infty}}\langle T(F),G \rangle{}_{\mathbb{D}^{\infty}}(=:\mathbb{E}[T(F)G])$ is the pairing or the generalized expectation of $T(F)\in \mathbb{D}^{-\infty}$ and $G\in \mathbb{D}^{\infty}$, and $\partial^{\alpha}T=\partial_{\alpha_1}\cdots \partial_{\alpha_k} T$ is understood as the distributional derivative sense.\\

We now discuss asymptotic expansion of Wiener functionals. 
For $\{ G_\varepsilon \}_{\varepsilon\in(0,1]} \subset \mathbb{D}^\infty$, we say $G_\varepsilon=O(\varepsilon^{r})$ in $\mathbb{D}^\infty$ if $\| G_\varepsilon \|_{k,p}=O(\varepsilon^{r})$ for all $k\in \mathbb{N}$ and $p>1$.  
Watanabe (1987) \cite{W87} shows that if a family of Wiener functionals $\{ F^{\varepsilon} \}_{\varepsilon \in (0,1]} \subset (\mathbb{D}^\infty)^e$ satisfies 
\begin{enumerate}
\item[(a)] 
\begin{eqnarray}
F^{\varepsilon,i} \sim F^{0,i}+\varepsilon F_1^{i}+\varepsilon^2 F_2^{i}+ \cdots + \ \ \mbox{in} \ \ \mathbb{D}^\infty, \ \ i=1,\ldots,e, \label{wat_condition_1}
\end{eqnarray}
where $F^0,F_1,F_2,\ldots \in (\mathbb{D}^\infty)^e$, in the sense that for any $m\geq 1$, 
\begin{eqnarray*}
F^{\varepsilon,i}-(F^{0,i}+\varepsilon F_1^{i}+\varepsilon^2 F_2^{i}+ \cdots + \varepsilon^m F_m^{i} )=O(\varepsilon^{m+1}) \ \ \mbox{in} \ \ \mathbb{D}^\infty, \ \ i=1,\ldots,e,
\end{eqnarray*}
\item[(b)] (the uniformly nondegenerate condition)
\begin{eqnarray}
\limsup_{\varepsilon \downarrow 0}\| \det(\sigma^{F^{\varepsilon}})^{-1} \|_{{\cal L}^p({\cal W})} <\infty \ \mbox{ for all} \ p>1, \label{wat_condition_2}
\end{eqnarray}
\end{enumerate}
then, for all $T \in {\cal S}'(\mathbb{R}^e)$, 
it holds that
\begin{eqnarray}
\mathbb{E}[T(F^{\varepsilon})]=a_0+\varepsilon a_1+\varepsilon^2 a_2+ \cdots \label{wat_expansion}
\end{eqnarray}
where
\begin{align}
&a_0=\mathbb{E}[T(F^{0})], \ \ a_1=\mathbb{E}[ \textstyle{\sum_{i=1}^e} \partial_i T(F^{0}) F_1^i ],  \nn\\ 
a_2=&\mathbb{E}[\textstyle{\sum_{i=1}^e} \partial_i T(F^{0}) F_2^i]+\mathbb{E}[ \textstyle{\frac{1}{2}\sum_{i_1,i_2=1}^e} \partial_{i_1}\partial_{i_2} T(F^{0}) F_1^{i_1}F_1^{i_2} ], \ \ \ldots. \label{wat_rep} 
\end{align}
\ \\

In this paper, we improve the conditions (\ref{wat_condition_1}), (\ref{wat_condition_2}) and the resulting expansion (\ref{wat_expansion}) 
with  the coefficients (\ref{wat_rep}) in Watanabe's expansion on the abstract Wiener space, which enables us to apply our asymptotic expansion in more general mathematical settings including solutions of rough differential equations and the functionals of fractional Brownian motions of irregular cases (i.e., the Hurst index $H<1/2$). 
We show a new {\it fractional order} expansion formula of $\mathbb{E}[f(F^{\varepsilon})]$ for a family of Wiener functionals $\{ F^{\varepsilon} \}_{\varepsilon \in (0,1]}$ in the following sense:
\begin{enumerate}
\item $F^{\varepsilon}$ has a fractional order expansion in $(\mathbb{D}^\infty)^e$ which is more general than (\ref{wat_condition_1}).
\item  It works under a weaker condition than the uniformly nondegenerate condition (\ref{wat_condition_2}).
\item  An asymptotic expansion is obtained as an extension of (\ref{wat_expansion}) with a new representation of expansion coefficients through iterative annihilation (Malliavin derivative) and creation (Skorohod integral) calculation with inner product of tensor products of the Hilbert (Cameron-Martin) space. 
Namely, the representation (\ref{wat_rep}) of the coefficients of Watanabe's expansion is generalized through a computation scheme with the Stroock-Taylor  formula, the chain rule of Malliavin derivative, the duality formula and the IBP formula, which can be applied to various problems.
\end{enumerate}
For the new expansion, we give the theoretical error including the uniform bound of $f$. \\ 

The first main result is as follows. 

\begin{thm}
\label{ae_thm}
Let $\{ F^{\varepsilon} \}_{\varepsilon \in (0,1]} \subset (\mathbb{D}^\infty)^e$ be a family of Wiener functionals such that $F^{\varepsilon}$ has an asymptotic expansion in $(\mathbb{D}^\infty)^e$: 
\begin{align}
F^{\varepsilon,i}\sim F^{0,i}+\varepsilon^{\kappa_1} F_1^{i}+\varepsilon^{\kappa_2} F_2^{i}+ \cdots  \ \ \mbox{in} \ \ \mathbb{D}^\infty, \ \ i=1,\ldots,e,
\end{align}
where $F^0,F_1,F_2,\ldots \in (\mathbb{D}^\infty)^e$ and $\{ \kappa_i ; i\in \mathbb{N}\}$ satisfies $0<\kappa_1<\kappa_2<\cdots$, in the sense that for any $m\geq 1$, 
\begin{align}
F^{\varepsilon,i}-(F^{0,i}+\varepsilon^{\kappa_1} F_1^{i}+\varepsilon^{\kappa_2} F_2^{i}+ \cdots + \varepsilon^{\kappa_m} F_m^{i} )=O(\varepsilon^{\kappa_{m+1}}) \ \ \mbox{in} \ \ \mathbb{D}^\infty, \ \ i=1,\ldots,e,
\end{align}
 and assume that the Malliavin covariance matrix $\sigma^{F^0}$ is invertible a.s. and 
\begin{align}
\| (\det \sigma^{F^0})^{-1} \|_{{\cal L}^p({\cal W})} <\infty, \label{weak_condition}
\end{align}
 for all $p >1$. Then, for $m\geq 1$, there exists $C>0$ such that 
\begin{align}
&\Big| \mathbb{E}[f(F^\varepsilon) ]- \Big\{ \mathbb{E}[f(F^0)]  \nonumber\\
& ~~~~ +\sum_{j=1}^m \varepsilon^{\nu_j} \sum_{k,\alpha,\beta,\gamma}^{(j)} \mathbb{E} \Big[ f(F^0) H_{\alpha \ast \gamma} \Big(F^0, \frac{1}{p!} \Big\langle DF^{0,\gamma_1} \otimes \cdots \otimes DF^{0,\gamma_p}, \mathbb{E}[D^p \prod_{i=1}^k F_{\beta_i}^{\alpha_i}] \Big\rangle_{{\cal H}^{\otimes p}} \Big) \Big] \Big\} \Big| \nonumber\\
&\leq C \| f \|_{\infty} \varepsilon^{\nu_{m+1}},    \label{main_expansion_formula}
\end{align} 
for any bounded measurable function $f:\mathbb{R}^e \to \mathbb{R}$  and $\varepsilon \in(0,1]$, where  
$\nu_\ell$, $\ell \in \mathbb{N}$ are all the elements of $\{ \textstyle{\sum_{i=1}^m} \kappa_{\beta_i}; \ \beta_1,\ldots,\beta_m \in \mathbb{N}, m \in \mathbb{N} \}$ in increasing order, and
\begin{align}
\sum_{k,\alpha,\beta,\gamma}^{(j)}=\sum_{ \substack{\beta=(\beta_1,\ldots,\beta_k) \in \mathbb{N}^k, k \in \mathbb{N},\\ \sum_{\ell=1}^k \kappa_{\beta_\ell}=\nu_j }} \ \sum_{\alpha=(\alpha_1,\ldots,\alpha_k)\in \{1,\ldots,e \}^k} \frac{1}{k!} \sum_{\gamma \in \{1,\ldots,e \}^{p},  p \geq 0 }.
\end{align}
Here, $\alpha \ast \gamma$ represents $\alpha \ast \gamma=(\alpha_1,\ldots,\alpha_k,\gamma_1,\ldots,\gamma_p)$ for $\alpha=(\alpha_1,\ldots,\alpha_k)$ and $\gamma=(\gamma_1,\ldots,\gamma_p)$, and we used the convention: if $p=0$, $\textstyle{\frac{1}{p!} \langle DF^{0,\gamma_1} \otimes \cdots \otimes DF^{0,\gamma_p}, \mathbb{E}[D^p G] \rangle_{{\cal H}^{\otimes p}}}=\mathbb{E}[G]$.
\end{thm}
\noindent\\
{\it Proof of Theorem \ref{ae_thm}}. 
We note that for all $\lambda \in [0,1]$, 
{{\begin{align*}
\Big|\det \sigma^{F^0+\lambda (F^{\varepsilon}-F^0)} - \det \sigma^{F^0} \Big|
 \leq (C \| D(F^{\varepsilon}-F^0) \|_{\cal H}^2 (\| DF^0 \|^2_{\cal H}+\| DF^{\varepsilon} \|^2_{\cal H})^{(2e-1)/2} )^{1/2} 
\end{align*}}}
for some $C>0$, and 
{{\begin{align*}
&\det \sigma^{F^0+\lambda (F^{\varepsilon}-F^0)} 
 \geq \det \sigma^{F^0} - (C \| D(F^{\varepsilon}-F^0) \|_{\cal H}^2 (\| DF^0 \|^2_{\cal H}+\| DF^{\varepsilon} \|^2_{\cal H})^{(2e-1)/2} )^{1/2}
\end{align*}}}
by (2.110) of \cite{BCR}. Let $\psi \in C_b^\infty(\mathbb{R})$, $0 \leq \psi \leq 1$ be given by 
\begin{align*}
\psi(x)={\bf 1}_{|x|\leq 1/8}+\exp \Big(1-(1/8)^2/((1/8)^2-(x-1/8)^2)\Big){\bf 1}_{1/8<|x|<1/4}, \ x \in \mathbb{R},
\end{align*}
and for $\varepsilon \in (0,1]$, let 
\begin{align*}
\eta^\varepsilon=\frac{C \| D(F^{\varepsilon}-F^0) \|_{\cal H}^2 (\| DF^0 \|^2_{\cal H}+\| DF^{\varepsilon} \|^2_{\cal H})^{(2e-1)/2}}{ (\det \sigma^{F^0})^2 }
\end{align*}
so that 
\begin{align}
\psi(\eta^\varepsilon) \neq 0 \ \ \ \mbox{implies} \ \ \ \det \sigma^{F^0+\lambda (F^{\varepsilon}-F^0)}\geq (1/2) \det \sigma^{F^0} \ \ \mbox{for all} \ \lambda \in [0,1]. \label{estimate_perturbed_mall_cov}
\end{align} 
For $k \in \mathbb{N}$ and $p > 1$, the Wiener functional $\eta^\varepsilon$ is bounded by  
$\big\| \eta^\varepsilon \big\|_{k,p} \leq C_1 \| F^{\varepsilon} - F^0 \|^2_{k+1,r} \leq C_2 \varepsilon^{2\kappa_1}$, 
for some $C_1,C_2,r>0$ depending on $k,p$, using the estimates: for $\ell \in \mathbb{N}$ and $q > 1$, there exist $C_3(q)>0,C_4(\ell,q)>0$ and $C_5(\ell,q)>0$ such that $\| (\det \sigma^{F^0})^{-1} \|_{{\cal L}^q({\cal W})} \leq C_3(q)$, $\| F^{\varepsilon} \|_{\ell,q} \leq C_4(\ell,q)$ and $\| F^0 \|_{\ell,q} \leq C_5(\ell,q)$. We have that for all $k \in \mathbb{N}$ and $p > 1$, there is $C>0$ such that $\| \psi(\eta^\varepsilon) \|_{k,p} \leq C$. By the properties of $\psi$, we can see that $1-\psi(\eta^\varepsilon) \neq 0$ implies $\eta^\varepsilon \geq 1/8$. Then we have  
\begin{align*}
\| 1-\psi(\eta^\varepsilon) \|_{{\cal L}^1({\cal W})} \leq \mu(\eta^\varepsilon \geq 1/8) \leq 2^{3r} \mathbb{E}[| \eta^\varepsilon |^r] \leq C(r) \varepsilon^{2\kappa_1 r}
\end{align*}
for arbitrary $r> 1$. Also, since for all $j \geq 1$, $\partial^j \psi(\eta^\varepsilon) \neq 0$ implies $1/8 < \eta^\varepsilon < 1/4$, we have for all $k \in \mathbb{N}$ and $p > 1$, 
\begin{align*}
\| 1-\psi(\eta^\varepsilon) \|_{k,p} = O(\varepsilon^{r})
\end{align*}
for arbitrary $r> 1$. 

Let $f \in {\cal S}(\mathbb{R}^e)$ be a bounded function. Consider the decomposition 
\begin{align}
\mathbb{E}[f(F^\varepsilon) ]=\mathbb{E}[f(F^\varepsilon) ( 1- {\psi(\eta^\varepsilon)} ) ]+\mathbb{E}[f(F^\varepsilon) \psi(\eta^\varepsilon) ]. \label{decomp}
\end{align}
 For the first term of the right-hand side of (\ref{decomp}), we have 
\begin{align*}
| \mathbb{E}[f(F^\varepsilon) ( 1-\psi(\eta^\varepsilon) ) ] | \leq \| f \|_{\infty} h(\varepsilon), 
\end{align*}
where $h(\varepsilon)=O(\varepsilon^{r})$ for any $r >0$. We next expand the second term of of the right-hand side of (\ref{decomp}). For $m \in \mathbb{N}$, let $N=N(m)\in \mathbb{N}$ such that $\kappa_1 (N+1)\geq \nu_{m+1}$. We have 
\begin{align}
\mathbb{E}[f(F^\varepsilon) \psi(\eta^\varepsilon) ] \nonumber
&=\mathbb{E}[f(F^0) \psi(\eta^\varepsilon) ]
+\sum_{i=1}^N \sum_{\alpha \in\{1,\ldots,e \}^i} \frac{1}{i!} \mathbb{E}[ \partial^{\alpha} f(F^0) \prod_{\ell=1}^i (F^{\varepsilon,\alpha_\ell}-F^{0,\alpha_\ell}) \psi(\eta^\varepsilon) ]+R^{\varepsilon}_{1,f} \nonumber\\
&=
\mathbb{E}[f(F^0)]+\sum_{j=1}^m \varepsilon^{\nu_j} \sum_{k,\alpha,\beta}^{(j)} \mathbb{E}[ \partial^{\alpha}f(F^0)\prod_{i=1}^kF_{\beta_i}^{\alpha_i} ]+R^{\varepsilon}_{1,f}+R^{\varepsilon}_{2,f} \label{trucation_expansion}
\end{align}
where $\textstyle{\sum_{k,\alpha,\beta}^{(j)}=\sum_{ \substack{\beta=(\beta_1,\ldots,\beta_k) \in \mathbb{N}, k \in \mathbb{N},\\ \sum_{\ell=1}^k \kappa_{\beta_\ell}=\nu_j }} \ \sum_{\alpha=(\alpha_1,\ldots,\alpha_k)\in \{1,\ldots,e \}^k} \frac{1}{k!}}$, $R^{\varepsilon}_{1,f}$ is given by 
\begin{align*}
R^{\varepsilon}_{1,f}&=\int_0^1 \frac{(1-\lambda)^{N}}{N!} \sum_{\alpha \in \{1,\ldots,e \}^{N+1}} \mathbb{E}[\partial^\alpha f(\widetilde{F^{\lambda,\varepsilon}}) \prod_{\ell=1}^{N+1} (F^{\varepsilon,\alpha_\ell}-F^{0,\alpha_\ell}) \psi(\eta^\varepsilon) ] d\lambda,
\end{align*}
with $\widetilde{F^{\lambda,\varepsilon}}=F^0+\lambda (F^\varepsilon-F^0)$, $\lambda \in [0,1]$, $\varepsilon \in(0,1]$, and $R^{\varepsilon}_{2,f}$ has the form: 
\begin{align*}
R^{\varepsilon}_{2,f}&= \sum_{\alpha \in \{1,\ldots,e \}^k,k\leq N} \mathbb{E}[ \partial^\alpha f(F^0) G_\alpha^\varepsilon \psi(\eta^\varepsilon) ]
+
\sum_{\alpha \in \{1,\ldots,e \}^k,k\leq N} \mathbb{E}[ \partial^\alpha f(F^0) \hat{G}_\alpha^\varepsilon (1-\psi(\eta^\varepsilon)) ]
\end{align*}
with $\{G_\alpha^\varepsilon\}_{\alpha \in \{1,\ldots,e \}^k,k\leq N,\varepsilon \in (0,1]},\{\hat{G}_\alpha^\varepsilon\}_{\alpha \in \{1,\ldots,e \}^k,k\leq N,\varepsilon \in (0,1]} \subset \mathbb{D}^\infty$ such that for all $k\leq N$ and $\alpha \in \{1,\ldots,e \}^k$, $G_\alpha^\varepsilon, \hat{G}_\alpha^\varepsilon$, $\varepsilon \in (0,1]$ satisfy for all $\ell \in \mathbb{N}$, $p > 1$, $\| G_\alpha^\varepsilon \|_{\ell,p}=O(\varepsilon^{\nu_{m+1}})$, $\|\hat{G}_\alpha^\varepsilon \|_{\ell,p}=O(\varepsilon^{|\alpha|\kappa_1})$.  
By (\ref{estimate_perturbed_mall_cov}), we have 
\begin{align*}
R^{\varepsilon}_{1,f}
=
\int_0^1 \frac{(1-\lambda)^{N}}{N!} \sum_{\alpha \in \{1,\ldots,e \}^{N+1}} \mathbb{E} \Big[f(\widetilde{F^{\lambda,\varepsilon}}) H_{\alpha}\Big( \widetilde{F^{\lambda,\varepsilon}}, \prod_{\ell=1}^{N+1} (F^{\varepsilon,\alpha_\ell}-F^{0,\alpha_\ell}) \psi(\eta^\varepsilon) \Big) \Big] d\lambda
\end{align*}
with the estimates: for $p\geq 1$, there exist $C>0$, $p_1,p_2,q>1$ and $k\in\mathbb{N}$ such that
\begin{align*}
\Big\| H_{\alpha}(\widetilde{F^{\lambda,\varepsilon}},\prod_{\ell=1}^{N+1} (F^{\varepsilon,\alpha_\ell}-F^{0,\alpha_\ell}) \psi(\eta^\varepsilon)) \Big\|_p \leq C \| (\det \sigma^{F^0})^{-1} \|_{p_1}^{q} \| \prod_{\ell=1}^{N+1} (F^{\varepsilon,\alpha_\ell}-F^{0,\alpha_\ell}) \psi(\eta^\varepsilon) \|_{k,p_2}
\end{align*}
by p102 of \cite{N}, and for $k\in\mathbb{N}$ and $p> 1$, $\| \textstyle{\prod_{\ell=1}^{N+1}} (F^{\varepsilon,\alpha_\ell}-F^{0,\alpha_\ell}) \|_{k,p}=O(\varepsilon^{\kappa_1 (N+1)})=O(\varepsilon^{\nu_{m+1}})$, $\| \psi(\eta^\varepsilon) \|_{k,p}=O(1)$. Then, there exists $C>0$ such that
\begin{align*}
|R^{\varepsilon}_{1,f}| \leq C \| f \|_{\infty} \varepsilon^{\nu_{m+1}},
\end{align*}
for all $\varepsilon \in (0,1]$. Also, we have the similar estimate for $R^{\varepsilon}_{2,f}$, i.e. there exists $C>0$ such that
\begin{align*}
|R^{\varepsilon}_{2,f}|&\leq \sum_{\alpha \in \{1,\ldots,e \}^k,k\leq N}  |\mathbb{E} [ f(F^0) H_{\alpha}(F^0, G_\alpha^\varepsilon \psi(\eta^\varepsilon)) ]| \\
&+\sum_{\alpha \in \{1,\ldots,e \}^k,k\leq N} |\mathbb{E} [ f(F^0) H_{\alpha}(F^0, \hat{G}_\alpha^\varepsilon (1-\psi(\eta^\varepsilon)) ]| \\
&\leq C\| f \|_{\infty}  \varepsilon^{\nu_{m+1}},
\end{align*}
for all $\varepsilon \in (0,1]$, since we have for $\alpha \in \{1,\ldots,e \}^k$, $k\leq N$, for all $\ell \in \mathbb{N}$, $p > 1$, $\| G_\alpha^\varepsilon \|_{\ell,p}=O(\varepsilon^{\nu_{m+1}})$, and for all $\ell \in \mathbb{N}$, $p > 1$, $\| 1- \psi(\eta^\varepsilon) \|_{\ell,p}=O(\varepsilon^q)$ for arbitrary $q>1$.

We give the representation of the expansion coefficients. While the similar computation in the error analysis can be applied to the expansion coefficients, we provide more useful representation for each coefficient for the practical computational purpose. 
Let $G=\textstyle{\prod_{i=1}^k} F_{\beta_i}^{\alpha_i}$. 
Then we have 
\begin{align} 
&\mathbb{E}[ \partial^{\alpha}f(F^0)G ] \nonumber\\
=&
\mathbb{E}[ \partial^{\alpha}f(F^0)] \mathbb{E}[ G ]+\mathbb{E}[ \partial^{\alpha}f(F^0) \sum_{p^\geq 1} \frac{1}{p!} \delta^p(\mathbb{E}[D^pG]) ] \label{eq1} \\
=&
\mathbb{E}[ \partial^{\alpha}f(F^0)] \mathbb{E}[ G ]+\sum_{p^\geq 1} \frac{1}{p!} \mathbb{E}[  \langle D^p \partial^{\alpha}f(F^0), \mathbb{E}[D^pG] \rangle_{{\cal H}^{\otimes p}} ] \label{eq2} \\
=&
\mathbb{E}[ \partial^{\alpha}f(F^0)] \mathbb{E}[ G ]+\sum_{p^\geq 1} \sum_{\gamma} \mathbb{E}[ \partial^{\alpha\ast \gamma}f(F^0) \frac{1}{p!} \langle DF^{0,\gamma_1}\otimes \cdots \otimes DF^{0,\gamma_p}, \mathbb{E}[D^pG] \rangle_{{\cal H}^{\otimes p}} ] \label{eq3} \\
=&
\sum_{p^\geq 0} \sum_{\gamma} \mathbb{E}[ f(F^0) H_{\alpha\ast \gamma}(F^0, \frac{1}{p!} \langle DF^{0,\gamma_1}\otimes \cdots \otimes DF^{0,\gamma_p}, \mathbb{E}[D^pG] \rangle_{{\cal H}^{\otimes p}} ) ],  \label{eq4}
\end{align}
where we applied the Stroock-Taylor formula (\ref{Stroock_formula}) in (\ref{eq1}), the duality formula (\ref{duality_formula}) in (\ref{eq2}), the chain rule of Malliavin derivative in (\ref{eq3}) and the IBP formula (\ref{IBP_formula}) in (\ref{eq4}). In the above, we used the convention: if $p=0$, $\textstyle{\frac{1}{p!} \langle DF^{0,\gamma_1} \otimes \cdots \otimes DF^{0,\gamma_p}, \mathbb{E}[D^p G] \rangle_{{\cal H}^{\otimes p}}}=\mathbb{E}[G]$.

Therefore, 
 for bounded measurable function $f:\mathbb{R}^e \to \mathbb{R}$, we have 
\begin{align*}
&\mathbb{E}[f(F^\varepsilon) ] = \mathbb{E}[f(F^0)]  \nonumber\\
& +\sum_{j=1}^m \varepsilon^{\nu_j} \sum_{k,\alpha,\beta,\gamma}^{(j)} \mathbb{E} \Big[ f(F^0) H_{\alpha \ast \gamma} \Big(F^0, \frac{1}{p!} \Big\langle DF^{0,\gamma_1} \otimes \cdots \otimes DF^{0,\gamma_p}, \mathbb{E}[D^p \prod_{i=1}^kF_{\beta_i}^{\alpha_i}] \Big\rangle_{{\cal H}^{\otimes p}} \Big) \Big] 
+\widetilde{R^{\varepsilon}_{f}},
\end{align*} 
with the estimate
\begin{align*}
|\widetilde{R^{\varepsilon}_{f}}| \leq C \| f \|_{\infty} \varepsilon^{\nu_{m+1}},
\end{align*}
for some $C>0$ independent of $f$ and $\varepsilon$.
 $\Box$\\

\begin{rem}\label{rem_section2}
The form of the factor $\langle DF^{0,\gamma_1} \otimes \cdots \otimes DF^{0,\gamma_p}, \mathbb{E}[D^p \textstyle{\prod_{i=1}^k} F_{\beta_i}^{\alpha_i}] \rangle_{{\cal H}^{\otimes p}}$ in the expansion coefficients in (\ref{main_expansion_formula}) in Theorem \ref{ae_thm} is crucial in applications as we only need the computation of Malliavin derivatives with an inner product computation on ${\cal H}^{\otimes p}$.  
In particular, it plays an important role in the derivation of the asymptotic expansion of the expectations of irregular functionals of solutions of rough differential equations in the next section. \\
\end{rem}

\begin{rem}\label{rem_section2}
Our condition (\ref{weak_condition}) weaker than  the uniformly nondegenerate condition (\ref{wat_condition_2}) is useful in various applications, since it can be easily checked without any complicated procedures or mathematical proofs. 
\end{rem}

\section{Asymptotic expansion formula of expectation of solution of rough differential equation driven by fractional Brownian motion with $H<1/2$}

In the section, we show asymptotic expansion formulas of expectations of a solution of a multidimensional rough differential equation driven by $d$-dimensional fractional Brownian motion with $H \in (1/3,1/2)$. Our setting mostly follows Cass and Lim (2019) \cite{CL}, Decreusefond and \"{U}st\"{u}nel (1999) \cite{DU}, Al\'os et al. (2001) \cite{AMN} and Nualart (2006) \cite{N}. 
The framework we consider in the section is a particular case of the Malliavin calculus developed in Section 2 on an abstract Wiener space.

Let ${\cal W}^d=C_0([0,T];\mathbb{R}^d):=\{ \omega:[0,T] \to \mathbb{R}^d; \ \omega \ \mbox{is continuous}, \ \omega(0)=0 \}$ with the supremum topology, and let $\mathscr{B}({\cal W}^d)$ be the Borel $\sigma$-field. 
Let $\mathbb{P}=\mu^d$ be the unique probability measure on $({\cal W}^d,\mathscr{B}({\cal W}^d))$ such that the canonical process $B^H=(B^{H,1},\ldots,B^{H,d})$ is a $d$-dimensional fractional Brownian motion with the Hurst index $H \in (1/3,1/2)$, that is, $\mathbb{E}[B_t^{H,i}B_s^{H,j}]=R(t,s) {\bf 1}_{i=j}$ for $t,s\in [0,T]$, $i,j=1,\ldots,d$ where
\begin{align}
R(t,s)=\frac{1}{2}\{ |s|^{2H}+|t|^{2H}-|t-s|^{2H} \}.
\end{align}
We denote by $\{ e_1,\ldots,e_d \}$ the canonical basis on $\mathbb{R}^d$. 
Let ${\cal H}^d$ be the Cameron-Martin space, the completion of the linear span of $\{ R(t, \cdot) e_j ; \ t\in [0,T], j=1,\ldots,d  \}$ with respect to the norm 
$\| \cdot \|_{{\cal H}^d}=\langle \cdot,\cdot \rangle_{{\cal H}^d}$ where 
\begin{align}
\langle R(t, \cdot) e_i,R(s, \cdot) e_j  \rangle_{{\cal H}^d}=R(t,s) {\bf 1}_{i=j}.
\end{align}
Let $K_{H}$ be the kernel for the singular case given by 
\begin{align}
K_{H}(t,s)=&
\Big(\frac{2H}{(1-2H)\beta(1-2H,H+1/2)}\Big)^{1/2} \Big[ \Big( \frac{t}{s} \Big)^{H-1/2}(t-s)^{H-1/2} \nn\\
& \hspace{8em} -(H-1/2) s^{1/2-H} \int_s^t u^{H-3/2}(u-s)^{H-1/2}du \Big]
\end{align}
where $\beta(\cdot,\cdot)$ denotes the Beta function. The Cameron-Martin space ${\cal H}^d$ is given by 
\begin{align}
{\cal H}^d=\{ f\in {\cal W}^d; \ \exists \tilde{f} \in L^2([0,T];\mathbb{R}^d) \ \mathrm{s.t.} \ f(t)=\int_0^t K_H(t,s) \tilde{f}(s) ds, \ t\geq 0  \}
\end{align}
(see Theorem 3.3.1 of Decreusefond and \"{U}st\"{u}nel (1999) \cite{DU} and Theorem 2.1 of Decreusefond (2001) \cite{D}). 
Let $\bar{\cal H}^d$ be the completion of the linear span of $\{ {\bf 1}_{[0,t]}(\cdot) e_j ; \ t\in [0,T], j=1,\ldots,d  \}$ with respect to the norm 
$\| \cdot \|_{\bar{\cal H}^d}=\langle \cdot,\cdot \rangle_{\bar{\cal H}^d}$ where 
\begin{align}
\langle {\bf 1}_{[0,t]}(\cdot) e_i, {\bf 1}_{[0,s]}(\cdot) e_j  \rangle_{\bar{\cal H}^d}=R(t,s) {\bf 1}_{i=j},
\end{align}
and then there exists an isomorphism $\Phi:\bar{\cal H}^d \to {\cal H}^d $ obtained from extending the map ${\bf 1}_{[0,t]}(\cdot)e_j \mapsto R(t,\cdot)e_j$, $t\in[0,T]$, $j=1,\ldots,d$ (see Definition 2.15 of Cass and Lim (2019) \cite{CL}).  
Thus, we have the following relationship: 
{\footnotesize{\begin{align}
\langle {\bf 1}_{[0,t]}e_i,{\bf 1}_{[0,s]}e_j \rangle_{\bar{\cal H}^d}=R(t,s) {\bf 1}_{i=j}
=\langle R(t,\cdot)e_i,R(s,\cdot)e_j \rangle_{{\cal H}^d}
=\langle K_H(t,\cdot) {\bf 1}_{[0,t]},K_H(s,\cdot) {\bf 1}_{[0,s]} \rangle_{L^2([0,T])} {\bf 1}_{i=j}.
\end{align}}}
On $({\cal W}^d,\bar{\cal H}^d,\mu^d)$, we can apply Malliavin calculus for fractional Brownian motion case. Recall that the map ${\bf 1}_{[0,t]}(\cdot) \mapsto B^H_t$ extends to a unique linear isometry $I$ from $\bar{\cal H}^d \to L_2({\cal W}^d)$ and $I(h)$ is a Gaussian random variable with mean $0$ and variance $\| h \|^2_{\bar{\cal H}^d}$. Let $\mathscr{S}_H({\cal W}^d)$ be the space of the functionals given by $\mathscr{S}_H({\cal W}^d)=\{F: {\cal W}^d \to \mathbb{R} \ ; \ F=f(I(h_1),\ldots,I(h_n))$, $n\in \mathbb{N}$, $f \in C_b^\infty(\mathbb{R}^n)$, $h_1,\ldots,h_n \in \bar{\cal H}^d\}$. For $F=f(I(h_1),\ldots,I(h_n)) \in \mathscr{S}({\cal W}^d)$, we define the Malliavin derivative $DF \in \bar{\cal H}^d$ as 
\begin{align}
DF=\sum_{i=1}^n (\partial_i f)(I(h_1),\ldots,I(h_n)) h_i,
\end{align}
or 
\begin{align}
D_{\ell,t} F=\sum_{i=1}^n (\partial_i f)(I(h_1),\ldots,I(h_n)) h_i^\ell(t), \ \ t\geq 0, \ \ \ell=1,\ldots,d.
\end{align}
The operator $D$ is a closable operator, and for $p> 1$, we define $\mathbb{D}_H^{1,p}=\overline{\mathscr{S}_H({\cal W}^d)}^{\| \cdot \|_{1,p}}$ where the norm $\| \cdot \|_{1,p}$ given by $\| F \|_{1,p}=\| F \|_{{\cal L}^p({\cal W}^d)}+\| DF \|_{{\cal L}^p({\cal W}^d;\bar{\cal H}^d)}$. Similarly, the higher-order Malliavin derivatives $D^k$ and the corresponding Sobolev spaces $\mathbb{D}_H^{k,p}$ can be defined iteratively. We define $\mathbb{D}_H^\infty=\cap_{k\in \mathbb{N},p> 1}\mathbb{D}_H^{k,p}$ and let $\mathbb{D}_H^{-\infty}$ be the dual space of $\mathbb{D}_H^\infty$.

Let $\mathrm{Dom}\delta_H^p=\{ u \in L^2({\cal W}^d; (\bar{\cal H}^d)^{\otimes p}); \ \exists C>0 \mbox{ s.t. } |\mathbb{E}[\langle D^pF,u\rangle_{(\bar{\cal H}^d)^{\otimes p}}]| \leq C \| F \|_{L^2({\cal W}^d)}, \forall F \in \mathbb{D}_H^{p,2} \}$. For $u=(u^1,\ldots,u^d) \in \mathrm{Dom}\delta_H(=\mathrm{Dom}\delta_H^1)$, there exists $\delta_H(u)=\textstyle{\sum_{i=1}^d} \delta_{H,i}(u^i) \in L^2({\cal W}^d)$ such that
\begin{align}
\mathbb{E}[\langle DF,u\rangle_{\bar{\cal H}^d}]=\mathbb{E}[F \delta_H(u) ].
\end{align}
On the setting, we still use notation on integration by parts (\ref{IBP_formula}) and generalized expectation in Section 2. 

Let ${\bf B}^H$ be the canonical geometric rough path lift and consider the following rough differential equation: 
\begin{eqnarray}
dX_t^{x}=V_0( X_t^{x} ) dt + V( X_t^{x} ) \circ d{\bf B}_t^{H}
\end{eqnarray}
starting from $X_0^{x}=x \in \mathbb{R}^e$, where $V_0 \in C_b^{\infty}(\mathbb{R}^e;\mathbb{R}^e)$ and $V=(V_1,\ldots,V_d) \in C_b^{\infty}(\mathbb{R}^e;\mathbb{R}^{e \times d})$ (see Friz and Victoir (2009) \cite{FV} or/and Friz and Hairer (2014) \cite{FH} for more details on rough differential equations). We assume the following elliptic condition. \\

\begin{assumption}
$V_1(x),\ldots,V_d(x)$ linearly span $\mathbb{R}^e$. \\
\end{assumption}

We introduce the scaling rough differential equation: 
\begin{eqnarray}
dX_t^{\varepsilon,x}=\varepsilon^{1/H}V_0( X_t^{\varepsilon,x} ) dt + \varepsilon V( X_t^{\varepsilon,x} ) \circ d{\bf B}_t^{H}, 
\ \ X_0^{\varepsilon,x}=x.
\end{eqnarray}
Note that $X_1^{t^H,x}$ and $X_t^{x}$ have the same probability law. We introduce $\| \alpha \|=\#\{ i ; \alpha_i \neq 0  \}+\#\{ i ; \alpha_i = 0  \} /H $ for a multi-index $\alpha=(\alpha_1,\ldots,\alpha_k) \in \{0,1,\ldots,d \}^k$, and define 
\begin{align}
V_i \varphi(\cdot)=\sum_{j=1}^e V_i^j(\cdot)\frac{\partial}{\partial x_j}\varphi(\cdot), \ \ i=0,1,\ldots,d,
\end{align}
 for a smooth function $\varphi: \mathbb{R}^e \to \mathbb{R}$. 
 
Let $\kappa_1<\kappa_2<\cdots$ be all the elements of 
\begin{align}
A=\{ k + \ell \times 1/H; \ k,\ell \in \mathbb{N} \cup \{ 0 \}, \ k+\ell\geq 1 \}
\end{align}
 in increasing order, that is, $\kappa_1=1$, $\kappa_2=2$, $\kappa_3=1/H$, $\kappa_4=3$, ... since $1/3<H<1/2$. 
 Here, the sequence $\kappa_1<\kappa_2<\cdots$ in $A$ is relevant to the following asymptotic expansion of $X_t^{\varepsilon,x}$:  
\begin{align}
X_t^{\varepsilon,x}=x+\varepsilon \sum_{i=1}^d V_i(x){B}_t^{H,i}+\sum_{k=2}^m \varepsilon^{\kappa_k} \sum_{\| \alpha \|=\kappa_k} V_{\alpha_1}\cdots V_{\alpha_{|\alpha|-1}}V_{\alpha_{|\alpha|}}(x) \mathbb{B}^{H}_{\alpha}(t)+R^\varepsilon_{m+1}(t,x)
\end{align}
whose expansion coefficients are obtained by (formal) Taylor expansion, which is justified as the stochastic Taylor expansion of Lyons-It\^o map (see Section 3.6 and Proposition 4.3 with (4.3) of Inahama (2016) \cite{I}), where $\mathbb{B}^{H}_{\alpha}$ is the iterated Stratonovich integral with respect to $B^H$ for a multi-index $\alpha$, i.e.  
\begin{align}
\mathbb{B}^{H}_{\alpha}(t)=\int_{0<t_1<\cdots <t_k<t}  \circ dB_{t_1}^{H,\alpha_1} \circ \cdots \circ dB_{t_k}^{H,\alpha_k}, \ \ \alpha=(\alpha_1,\ldots,\alpha_k)\in \{0,1,\ldots,d\}^k
\end{align} 
and $R^{\varepsilon,j}_{m+1}(t,x)=(R^{\varepsilon,j}_{m+1}(t,x),\ldots,R^{\varepsilon,e}_{m+1}(t,x))$ is the residual satisfying  $R^{\varepsilon,j}_{m+1}(t,x)=O(\varepsilon^{m+1})$ in $\mathbb{D}_H^\infty$ for $j=1,\ldots,e$. By the equation (7) (or Theorem 6.1 and Theorem 6.3) of Cass and Lim (2019) \cite{CL} and the equation (4) in Song and Tindel (2022) \cite{ST}, we are able to transform Stratonovich integrals (in rough path sense) into Skorohod integrals in our setting, and thus the all expansion terms are Malliavin differentiable at all orders. 

Let 
\begin{align}
F_t^\varepsilon :=& (X^{\varepsilon,x}_t-x)/\varepsilon \\
=&\sum_{i=1}^d V_i(x){B}_t^{H,i}+\sum_{k=2}^m \varepsilon^{\kappa_k-1} \sum_{\| \alpha \|=\kappa_k} V_{\alpha_1}\cdots V_{\alpha_{|\alpha|-1}}V_{\alpha_{|\alpha|}}(x) \mathbb{B}^{H}_{\alpha}(t)+ R^\varepsilon_{m+1}(t,x).
\end{align} 
Define $F_t^0:=\sum_{i=1}^d V_i(x){B}_t^{H,i}$ and 
\begin{align}
F_{t,k}^{\ell}:=\sum_{\| \alpha \|=\kappa_k} V_{\alpha_1}\cdots V_{\alpha_{{|\alpha|}-1}}V_{\alpha_{|\alpha|}}^{\ell}(x) \mathbb{B}^{H}_{\alpha}(t), \ \ t\geq 0, \ \ell=1,\ldots,e, \ k \in \mathbb{N}.
\end{align}
Under {\bf Assumption 1}, we easily check that for all $p\geq 1$, there exists $C>0$ such that
\begin{align}
\| (\det \sigma^{F_t^0})^{-1}\|_{{\cal L}^p({\cal W}^d)} \leq \frac{C}{t^{2He}}. \label{det_mall_2}
\end{align}

Let 
\begin{align}
\bar{X}_t^{x,\varepsilon}:=x + \varepsilon F_t^0=x+\varepsilon \sum_{i=1}^d V_i(x){B}_t^{H,i}, \ \ t>0, \ \varepsilon\in (0,1].
\end{align}
 
We have the following new expansion as an application of Theorem \ref{ae_thm} by taking $\bar{X}_t^{x,\varepsilon}|_{t=1,\varepsilon=t^H}$, whose expansion coefficients are more simplified (see Remark \ref{rem_general_expansion_fbm} below) than those in Theorem \ref{ae_thm} under the setting. 

\begin{thm}\label{general_expansion_formula_rde}
For $m\geq 1$, there exists $C>0$ such that 
\begin{align}
&\Big| \mathbb{E}[f(X_t^x) ]- \Big\{ \mathbb{E}[f(\bar{X}_1^{x,t^H})]  +\sum_{j=1}^m t^{H\nu_j} \sum_{k,\alpha,\beta,\gamma}^{(j)'} \nonumber\\
& ~~~~~~~~~~~~~
 \mathbb{E} \Big[ f(\bar{X}_1^{x,t^H}) H_{\alpha}\Big(F^0_1, H_{\gamma} \Big(B_1^H, \frac{1}{p!} \Big\langle DB_1^{H,\gamma_1} \otimes \cdots \otimes DB_1^{H,\gamma_p}, \mathbb{E}[D^p \prod_{i=1}^kF_{1,\beta_i}^{\alpha_i}] \Big\rangle_{(\bar{\cal H}^d)^{\otimes p}} \Big)\Big) \Big] \Big\} \Big| \nonumber\\
&\leq C \| f \|_{\infty} t^{H\nu_{m+1}}, \label{ste_1}
\end{align} 
for any bounded measurable function $f:\mathbb{R}^e \to \mathbb{R}$ and $t \in(0,1]$, where $\nu_\ell$, $\ell \in \mathbb{N}$ are all the elements of $\bar{A}=\{ \textstyle{\sum_{\ell=1}^m} (\kappa_{\beta_\ell}-1) ; \ \beta_1,\ldots,\beta_m \geq 2, \ m \in \mathbb{N} \}$ in increasing order, and 
\begin{align}
\sum_{k,\alpha,\beta,\gamma}^{(j)'}=\sum_{ \substack{\beta=(\beta_1,\ldots,\beta_k)\in (\mathbb{N}\setminus \{ 1\})^k, k \in \mathbb{N},\\ \sum_{\ell=1}^k (\kappa_{\beta_\ell}-1)=\nu_j }} \ \sum_{\alpha=(\alpha_1,\ldots,\alpha_k)\in \{1,\ldots,e \}^k} \frac{1}{k!} \sum_{\gamma \in \{1,\ldots,d \}^{p},  p \geq 0 }.
\end{align}
\end{thm}
\ \\
\begin{rem}
Note that the index set $\bar{A}$ in Theorem \ref{general_expansion_formula_rde} is relevant to the asymptotic expansion of $\mathbb{E}[\delta_y(F_t^\varepsilon)]=\mathbb{E}[\delta_y((X^{\varepsilon,x}_t-x)/\varepsilon)]$. \\
\end{rem}

\noindent
{\it Proof of Theorem \ref{general_expansion_formula_rde}}. \ 
Let $f \in {\cal S}(\mathbb{R}^e)$ be a bounded function. Note that we have
\begin{align}
\mathbb{E}[ f(X_t^x) ]=\int_{\mathbb{R}^e} f(x+\varepsilon y) {}_{\mathbb{D}_H^{-\infty}} \langle \delta_y(F_1^\varepsilon),1 \rangle_{\mathbb{D}_H^{\infty}} dy \Big|_{\varepsilon=t^H}. \label{watanabe_representation} 
\end{align}
By applying the proof of Theorem \ref{ae_thm}, it holds that
\begin{align}
& {}_{\mathbb{D}_H^{-\infty}} \langle \delta_y(F_1^\varepsilon),1 \rangle  {}_{\mathbb{D}_H^{\infty}} = 
{}_{\mathbb{D}_H^{-\infty}} \langle \delta_y (F_1^\varepsilon), {\psi(\eta^\varepsilon)}  \rangle_{\mathbb{D}_H^{\infty}}+{}_{\mathbb{D}_H^{-\infty}} \langle \delta_y (F_1^\varepsilon), ( 1- {\psi(\eta^\varepsilon)} ) \rangle_{\mathbb{D}_H^{\infty}} \nonumber\\
= & {}_{\mathbb{D}_H^{-\infty}} \langle \delta_y(F_1^0),1 \rangle {}_{\mathbb{D}_H^{\infty}} +\sum_{j=1}^m \varepsilon^{\nu_j} \sum_{k,\alpha,\beta}^{(j)} {}_{\mathbb{D}_H^{-\infty}} \Big\langle \partial^{\alpha} \delta_y(F_1^0), \prod_{i=1}^kF_{1,\beta_i}^{\alpha_i} \Big\rangle_{\mathbb{D}_H^{\infty}} +R_{0,\delta_y}(\varepsilon)+R_{1,\delta_y}(\varepsilon)+R_{2,\delta_y}(\varepsilon), \label{delta_expansion}
\end{align} 
where $\textstyle{\sum_{k,\alpha,\beta}^{(j)}=\sum_{ \substack{\beta=(\beta_1,\ldots,\beta_k)\in (\mathbb{N}\setminus \{ 1 \})^k, k \in \mathbb{N},\\ \sum_{\ell=1}^k (\kappa_{\beta_\ell}-1)=\nu_j }} \ \sum_{\alpha=(\alpha_1,\ldots,\alpha_k)\in \{1,\ldots,e \}^k} \frac{1}{k!}}$, 
\begin{align*}
\psi(x)=&\textstyle{{\bf 1}_{|x|\leq 1/8}+\exp (1-(1/8)^2/((1/8)^2-(x-1/8)^2)){\bf 1}_{1/8<|x|<1/4}}, \ \ x \in \mathbb{R},
\end{align*} 
\begin{align*}
\eta^\varepsilon=&\frac{C \| D(F_1^{\varepsilon}-F_1^0) \|_{\bar{\cal H}^d}^2 (\| DF_1^0 \|^2_{\bar{\cal H}^d}+\| DF_1^{\varepsilon} \|^2_{\bar{\cal H}^d})^{(2e-1)/2}}{ (\det \sigma^{F_1^0})^2 }, \ \ \varepsilon\in (0,1],
\end{align*} 
\begin{align*}
R_{0,\delta_y}(\varepsilon)={}_{\mathbb{D}_H^{-\infty}} \langle \delta_y (F_1^\varepsilon), ( 1- {\psi(\eta^\varepsilon)} ) \rangle_{\mathbb{D}_H^{\infty}},
\end{align*}
\begin{align*}
R_{1,\delta_y}(\varepsilon)=\int_0^1 \frac{(1-\lambda)^{N}}{N!} \sum_{\alpha \in \{1,\ldots,e \}^{N+1}} {}_{\mathbb{D}_H^{-\infty}} \Big\langle \delta_y(\widetilde{F_1^{\lambda,\varepsilon}}), H_{\alpha}\Big( \widetilde{F_1^{\lambda,\varepsilon}}, \prod_{\ell=1}^{N+1} (F_1^{\varepsilon,\alpha_\ell}-F_1^{0,\alpha_\ell}) \psi(\eta^\varepsilon) \Big) \Big\rangle_{\mathbb{D}_H^{\infty}} d\lambda
\end{align*}
with a natural number $N$ such that $\kappa_1 (N+1)\geq \nu_{m+1}$ and $\textstyle{\widetilde{F_1^{\lambda,\varepsilon}}=F_1^0+\lambda (F_1^\varepsilon-F_1^0)}$, $\lambda \in [0,1]$, $\varepsilon \in(0,1]$, 
and 
\begin{align*}
R_{2,\delta_y}(\varepsilon)=&
\sum_{\alpha \in \{1,\ldots,e \}^k,k\leq N} {}_{\mathbb{D}_H^{-\infty}} \langle \delta_y(F_1^0), H_{\alpha}(F_1^0, G_\alpha^\varepsilon \psi(\eta^\varepsilon) ) \rangle_{\mathbb{D}_H^{\infty}}\\
&+
\sum_{\alpha \in \{1,\ldots,e \}^k,k\leq N} {}_{\mathbb{D}_H^{-\infty}} \langle \delta_y(F_1^0), H_{\alpha}(F_1^0, \hat{G}_\alpha^\varepsilon (1-\psi(\eta^\varepsilon)) ) \rangle_{\mathbb{D}_H^{\infty}}
\end{align*}
with $\{G_\alpha^\varepsilon\}_{\alpha \in \{1,\ldots,e \}^k,k\leq N,\varepsilon \in (0,1]},\{ \hat{G}_\alpha^\varepsilon\}_{\alpha \in \{1,\ldots,e \}^k,k\leq N,\varepsilon \in (0,1]} \subset \mathbb{D}_H^\infty$ such that for any $k\leq N$ and multi-index $\alpha \in \{1,\ldots,e \}^k$, $G_\alpha^\varepsilon$, $\hat{G}_\alpha^\varepsilon$, $\varepsilon \in (0,1]$ satisfy for all $\ell \in \mathbb{N}$, $p > 1$, $\| G_\alpha^\varepsilon \|_{\ell,p}=O(\varepsilon^{\nu_{m+1}})$, $\|\hat{G}_\alpha^\varepsilon \|_{\ell,p}=O(\varepsilon^{|\alpha|\kappa_1})$. Note that we have for $k\in \mathbb{N}$, $p > 1$,
\begin{align*}
\| 1-\psi(\eta^\varepsilon) \|_{k,p} = O(\varepsilon^r)
\end{align*}
for arbitrary $r > 1$, and for all $k\in\mathbb{N}$, $p> 1$, $\| \psi(\eta^\varepsilon) \|_{k,p}=O(1)$. 
By (\ref{watanabe_representation}) and (\ref{delta_expansion}), it holds that 
\begin{align}
\mathbb{E}[f(X_t^x) ]=& \mathbb{E}[f(\bar{X}_1^{x,t^H})] +\sum_{j=1}^m t^{H\nu_j} \sum_{k,\alpha,\beta}^{(j)} \mathbb{E} \Big[ f(\bar{X}_1^{x,t^H}) H_{\alpha}\Big(F_1^0, \prod_{i=1}^kF_{1,\beta_i}^{\alpha_i} \Big) \Big]  \nonumber\\
&\hspace{15em}+\widetilde{R_{0,f}} (t^H)+\widetilde{R_{1,f}} (t^H)+\widetilde{R_{2,f}} (t^H), \label{E1}
\end{align} 
where $\widetilde{R_{0,f}}(\varepsilon)=\mathbb{E} [ f(x+\varepsilon F_1^\varepsilon) ( 1- {\psi(\eta^\varepsilon)} ) ]$, 
\begin{align*}
\widetilde{R_{1,f}}(\varepsilon)=\int_0^1 \frac{(1-\lambda)^{N}}{N!} \sum_{\alpha \in \{1,\ldots,e \}^{N+1}} \mathbb{E} \Big[ f(x+ \varepsilon \widetilde{F_1^{\lambda,\varepsilon}}) H_{\alpha}\Big( \widetilde{F_1^{\lambda,\varepsilon}}, \prod_{\ell=1}^{N+1} (F_1^{\varepsilon,\alpha_\ell}-F_1^{0,\alpha_\ell}) \psi(\eta^\varepsilon) \Big) \Big] d\lambda
\end{align*}
and
\begin{align*}
\widetilde{R_{2,f}}(\varepsilon)=&
\sum_{\alpha \in \{1,\ldots,e \}^k,k\leq N} \mathbb{E} [ f(x+\varepsilon F_1^0) H_{\alpha}(F_1^0, G_\alpha^\varepsilon \psi(\eta^\varepsilon) ) ]\\
&+
\sum_{\alpha \in \{1,\ldots,e \}^k,k\leq N} \mathbb{E} [ f(x+\varepsilon F_1^0) H_{\alpha}(F_1^0, \hat{G}_\alpha^\varepsilon (1-\psi(\eta^\varepsilon)) ) ]
\end{align*}
satisfy that there exists $C>0$ independent of $f$ such that 
\begin{align}
\Big|\sum_{k=0}^2 \widetilde{R_{k,f}}(\varepsilon) \Big| \leq C \| f \|_{\infty} \varepsilon^{\nu_{m+1}}, \label{error_ep_and_t}
\end{align}
for all $\varepsilon \in (0,1]$.  Note that one has
\begin{align}
\mathbb{E}[ f(\bar{X}_1^{x,\varepsilon}) H_\alpha(F_1^0,G) ]=&\mathbb{E}[ f(\bar{X}_1^{x,\varepsilon}) \varepsilon^{|\alpha|} H_\alpha(\bar{X}_1^{x,\varepsilon},G) ]= \varepsilon^{|\alpha| }  \mathbb{E}[ \partial^\alpha f(\bar{X}_1^{x,\varepsilon}) G ] \nonumber\\
=& \varepsilon^{|\alpha| } \int_{\mathbb{R}^e} \partial^\alpha f(x+\varepsilon V(x)y) \ {}_{\mathbb{D}_H^{-\infty}} \langle \delta_y(B_1^{H}),G \rangle {}_{\mathbb{D}_H^{\infty}}  dy,
\end{align}
for any $f \in {\cal S}(\mathbb{R}^e)$, $\varepsilon \in (0,1]$, $G \in \mathbb{D}_H^\infty$ and multi-index $\alpha$.
Since we have 
\begin{align} 
&{}_{\mathbb{D}_H^{-\infty}} \langle \delta_y (B_1^{H}),G \rangle_{\mathbb{D}_H^{\infty}} \nonumber \\
=&
 \sum_{\gamma \in \{1,\ldots,d \}^p, p^\geq 0} {}_{\mathbb{D}_H^{-\infty}} \Big\langle \delta_y(B_1^{H}),H_{\gamma}\Big(B_1^{H}, \frac{1}{p!} \Big\langle DB_1^{H,\gamma_1}\otimes \cdots \otimes DB_1^{H,\gamma_p}, \mathbb{E}[D^pG] \Big\rangle_{(\bar{\cal H}^d)^{\otimes p}} \Big) \Big\rangle {}_{\mathbb{D}_H^{\infty}} \label{brownian_rep}
\end{align}
by the similar argument in (\ref{eq1})--(\ref{eq4}), it holds that 
\begin{align}
&\mathbb{E}[ f(\bar{X}_1^{x,\varepsilon}) H_\alpha(F_1^0,G) ]\nonumber\\
=& \varepsilon^{|\alpha| }\sum_{\gamma \in \{1,\ldots,d \}^p, p^\geq 0}  \mathbb{E} \Big[ \partial^\alpha f(\bar{X}_1^{x,\varepsilon}) H_{\gamma} \Big( B_1^{H}, \frac{1}{p!} \Big\langle DB_1^{H,\gamma_1}\otimes \cdots \otimes DB_1^{H,\gamma_p}, \mathbb{E}[D^pG] \Big\rangle_{(\bar{\cal H}^d)^{\otimes p}} \Big) \Big]\nonumber\\
=& \sum_{\gamma \in \{1,\ldots,d \}^p, p^\geq 0}  \mathbb{E} \Big[ f(\bar{X}_1^{x,\varepsilon}) H_\alpha \Big(F_1^0,H_{\gamma} \Big(B_1^{H}, \frac{1}{p!} \Big\langle DB_1^{H,\gamma_1}\otimes \cdots \otimes DB_1^{H,\gamma_p}, \mathbb{E}[D^pG] \Big\rangle_{(\bar{\cal H}^d)^{\otimes p}} \Big) \Big) \Big]. \label{rep1}
\end{align}
Then, by (\ref{E1}), (\ref{rep1}) and the error estimate (\ref{error_ep_and_t}), 
 we have
\begin{align}
&\Big| \mathbb{E}[f(X_t^x) ]- \Big\{ \mathbb{E}[f(\bar{X}_1^{x,t^H})]  +\sum_{j=1}^m t^{H\nu_j} \sum_{k,\alpha,\beta,\gamma}^{(j)'} \nonumber\\
& ~~~~~~~~~~~~~
 \mathbb{E} \Big[ f(\bar{X}_1^{x,t^H}) H_{\alpha}\Big(F^0_1, H_{\gamma} \Big(B_1^H, \frac{1}{p!} \Big\langle DB_1^{H,\gamma_1} \otimes \cdots \otimes DB_1^{H,\gamma_p}, \mathbb{E}[D^p \prod_{i=1}^kF_{1,\beta_i}^{\alpha_i}] \Big\rangle_{(\bar{\cal H}^d)^{\otimes p}} \Big)\Big) \Big] \Big\} \Big| \nonumber\\
&\leq C \| f \|_{\infty} t^{H\nu_{m+1}},
\end{align} 
for any bounded measurable function $f:\mathbb{R}^e \to \mathbb{R}$ and $t \in(0,1]$. $\Box$ \\

\begin{rem}\label{rem_general_expansion_fbm}
The weight 
\begin{align}
H_{\alpha}\Big(F^0_1,H_{\gamma} \Big(B_1^H, \frac{1}{p!} \Big\langle DB_1^{H,\gamma_1} \otimes \cdots \otimes DB_1^{H,\gamma_p}, \mathbb{E}[D^p \prod_{i=1}^kF_{1,\beta_i}^{\alpha_i}] \Big\rangle_{(\bar{\cal H}^d)^{\otimes p}} \Big)\Big) 
\end{align}
 in (\ref{ste_1}) in Theorem \ref{general_expansion_formula_rde} is simpler than the weight: $$H_{\alpha \ast \gamma}\Big(F^0_1,\frac{1}{p!} \Big\langle DF_1^{0,\gamma_1} \otimes \cdots \otimes DF_1^{0,\gamma_p}, \mathbb{E}[D^p \prod_{i=1}^kF_{1,\beta_i}^{\alpha_i}] \Big\rangle_{(\bar{\cal H}^d)^{\otimes p}} \Big).$$ This is due to
the difference between the forms of Malliavin derivatives and the Malliavin covariance matrices of $F^0_1$ and $B_1^H$, respectively, i.e. $$\sigma^{F^0_1} = [ \textstyle{\sum_{k=1}^d} V_k^i(x) V_k^j(x)]_{1\leq i,j\leq e} \ \mbox{and} \ \sigma^{B_1^H}=[\delta_{i,j} ]_{1\leq i,j\leq d}.$$ Through the asymptotic expansion formulas given in Theorem \ref{general_expansion_formula_rde}, we can reduce the $|\gamma|$-times inverse Malliavin covariance matrix computation of $F_1^0$ in $$H_{\alpha \ast \gamma} \Big(F_1^0,G \Big) = H_{(({\alpha \ast \gamma})_k)} \Big( F_1^0,H_{({(\alpha \ast \gamma)}_1,\ldots,{(\alpha \ast \gamma)}_{k-1})} \Big(F_1^0,G \Big) \Big)$$ with $$H_{(i)} \Big(F_1^0,G \Big)=\sum_{j=1}^e \delta_H \Big( (\sigma^{F_1^0})^{-1}_{ij} D F_1^0 G \Big), \  i=1,\ldots,e.$$ We will see the effect in Theorem \ref{expansion_rde} below. \\
\end{rem}
 
Theorem \ref{general_expansion_formula_rde} enables us to give more explicit form of the expansion in each specific order of approximation without using complicated fractional calculus, which cannot be obtained by the previous approaches in the literature.
We only need an inner product computation on $(\bar{\cal H}^d)^{\otimes p}$ with IBP formula after we compute the Malliavin derivatives of $\textstyle{\prod_{i=1}^k}F_{1,\beta_i}^{\alpha_i}$ in the derivation of the asymptotic expansion in Theorem \ref{general_expansion_formula_rde}. 
As a consequence, all expansion terms are obtained as polynomials of fractional Brownian motion for multidimensional system of rough differential equations.\\

We have the following concrete asymptotic expansions as a main result of the paper. 
\begin{thm}\label{expansion_rde}
We have 
\begin{eqnarray}
&&\mathbb{P} ( X_t^x \leq y) = \mathbb{P}(\bar{X}_1^{x,t^H} \leq y) \nn \\
&& \hspace{1.5em}
 + t^{H} \mathbb{E} \Bigl[ {\bf 1}_{ \left\{ \bar{X}_1^{x,t^H} \leq y \right\} }
 \sum_{j_1,j_2=1}^e \sum_{i_1,i_2,i_3=1}^d V_{i_1}V_{i_2}^{j_1}(x) V_{i_3}^{j_2}(x) A^{-1}_{j_1,j_2}(x) \nn \\
&& \hspace{4.5em} \frac{1}{2} \{  B_1^{H,i_1}B_1^{H,i_2}B_1^{H,i_3} - B_1^{H,i_1} {\bf 1}_{i_2=i_3\neq 0} - B_1^{H,i_2}  {\bf 1}_{i_1=i_3\neq 0}  \} \Bigr]  + O(t^{1-H}) \ \ 
\end{eqnarray}
where $A(x)=\sum_{i=1}^d V_i(x) \otimes V_i(x) $. Moreover, we have 
\begin{eqnarray}
&&\mathbb{P}( X_t^x \leq y) = \mathbb{P}(\bar{X}_1^{x,t^H} \leq y) \nn \\
&& \hspace{1.5em}
 + t^{H} \mathbb{E} \Bigl[ {\bf 1}_{ \left\{ \bar{X}_1^{x,t^H} \leq y \right\} }
 \sum_{j_1,j_2=1}^e \sum_{i_1,i_2,i_3=1}^d V_{i_1}V_{i_2}^{j_1}(x) V_{i_3}^{j_2}(x) A^{-1}_{j_1,j_2}(x) \nn \\
&& \hspace{4.5em} \frac{1}{2} \{  B_1^{H,i_1}B_1^{H,i_2}B_1^{H,i_3} - B_1^{H,i_1} {\bf 1}_{i_2=i_3\neq 0} - B_1^{H,i_2}  {\bf 1}_{i_1=i_3\neq 0}  \} \Bigr] \nn \\
&& \hspace{1.5em} + t^{1-H} \mathbb{E} \Bigl[ {\bf 1}_{ \left\{ \bar{X}_1^{x,t^H} \leq y \right\} }
 \sum_{j_1,j_2=1}^e \sum_{i_1=1}^d V_{0}^{j_1}(x) V_{i_1}^{j_2}(x) A^{-1}_{j_1,j_2}(x) B_1^{H,i_1} \Bigr]  \ +  O(t^{2H}). 
\end{eqnarray}
\end{thm}
\noindent\\

\begin{rem}
The weight in Theorem \ref{expansion_rde} is given by the polynomial of Brownian motions but do not have Hermite polynomial structure due to the property of the geometric rough path integral. 
\end{rem}
\begin{rem}
The expansion in Theorem \ref{expansion_rde} is implemented by a simple numerical scheme. We will see it in the end of this section.   
\end{rem}
\begin{rem}
When $e=d$, the formula is obtained in a simple way using the inverse matrix of $V(\cdot)$.   
\end{rem}

\noindent
{\it Proof of Theorem \ref{expansion_rde}}. 
Hereafter, we use a notation $\int_0^t u_s^i \delta B_s^{H,i}:=\delta_{H,i}(u^i {\bf 1}_{[0,t]}(\cdot))$, $i=1,\ldots,d$ for $u=(u^1,\ldots,u^d) \in \mathrm{Dom}\delta_H$. In order to obtain more explicit form of the $O(t^{1-H})$-expansion, we trace the derivation in Theorem \ref{general_expansion_formula_rde} and compute the following term: 
\begin{eqnarray*}
&& t^H \mathbb{E}[  f ( \bar{X}_1^{x,t^H} ) H_{(j)}(F_1^0,V_{i_1}V_{i_2}^{j}(x) \mathbb{B}^{H}_{(i_1,i_2)}(1) ) ]\\
&=&
t^H \mathbb{E}[  f ( \bar{X}_1^{x,t^H} ) t^H H_{(j)}(\bar{X}_1^{x,t^H},V_{i_1}V_{i_2}^{j}(x) \mathbb{B}^{H}_{(i_1,i_2)}(1) ) ]\\
&=& t^{2H} \mathbb{E}[ \partial_{j} f ( \bar{X}_1^{x,t^H} ) V_{i_1}V_{i_2}^{j}(x) \mathbb{B}^{H}_{(i_1,i_2)}(1) ]\\
&=& t^{2H} \int_{\mathbb{R}^d} \partial_{j} f ( x+t^H V(x)y ) V_{i_1}V_{i_2}^{j}(x) \mathbb{E}[\delta_y(B_1^{H}) \mathbb{B}^{H}_{(i_1,i_2)}(1) ] dy, 
\end{eqnarray*}
for $j=1,\ldots,e$ and $i_1,i_2=1,\ldots,d$. We analyze $\mathbb{E}[\delta_y(B_t^{H}) \mathbb{B}^{H}_{(i_1,i_2)}(t) ]$ for $t>0$. By the equation (7) (or Theorem 6.1 and Theorem 6.3) of Cass and Lim (2019) \cite{CL} and the equation (4) in Song and Tindel (2021) \cite{ST}, it holds that
\begin{eqnarray}
\mathbb{B}^{H}_{(i_1,i_2)}(t)=\int_0^t \int_0^{t_{2}} \circ dB_{t_1}^{H,i_1} \circ dB_{t_2}^{H,i_2}&=&\int_0^t B_{s}^{H,i_1} \delta B_s^{H,i_2}+H \int_0^t s^{2H-1}ds {\bf 1}_{i_1=i_2} \nn\\
&=&\int_0^t B_{s}^{H,i_1} \delta B_s^{H,i_2}+\frac{1}{2} t^{2H} {\bf 1}_{i_1=i_2}.
\end{eqnarray}
Then we have
\begin{eqnarray}
\mathbb{E}[\delta_y(B_t^{H}) \mathbb{B}^{H}_{(i_1,i_2)}(t) ]=
\mathbb{E}[\delta_y(B_t^{H}) \int_0^t B_{s}^{H,i_1} \delta B_s^{H,i_2}] + \mathbb{E} [\delta_y(B_t^{H}) ] \frac{1}{2} t^{2H}{\bf 1}_{i_1=i_2}. \label{exp_iterated_integral}
\end{eqnarray}
We compute the Malliavin derivatives of $\int_0^t B_{s}^{H,i_1} \delta B_s^{H,i_2}$ to obtain each term constituting the Stroock-Taylor formula. We note that
\begin{align}
\mathbb{E}[ D^0 \int_0^t B_{s}^{H,i_1} \delta B_s^{H,i_2}]=\mathbb{E}[ \int_0^t B_{s}^{H,i_1} \delta B_s^{H,i_2}]=0.
\end{align}
First we compute $\mathbb{E}[D \int_0^t B_{s}^{H,i_1} \delta B_s^{H,i_2}]=(\mathbb{E}[D_{1,\cdot}\int_0^t B_{s}^{H,i_1} \delta B_s^{H,i_2}],\ldots,\mathbb{E}[D_{d,\cdot}\int_0^t B_{s}^{H,i_1} \delta B_s^{H,i_2}])$. For $\ell=1,\ldots,d$, we have 
\begin{eqnarray}
D_{\ell,r}\int_0^t B_{s}^{H,i_1} \delta B_s^{H,i_2}&=&B_{r}^{H,i_1} {\bf 1}_{\ell=i_2}+ \int_{r}^t D_{\ell,r}B_{s}^{H,i_1} \delta B_s^{H,i_2} \nn\\
&=&B_{r}^{H,i_1} {\bf 1}_{\ell=i_2} + (B_{t}^{H,i_2}-B_{r}^{H,i_2}) {\bf 1}_{\ell=i_1} \ \ \ \mbox{for} \ r \leq t \ \ 
 \label{first_mall_deriv}
\end{eqnarray}
and then 
\begin{align}
\mathbb{E}[D_{\ell,r}\int_0^t B_{s}^{H,i_1} \delta B_s^{H,i_2}]=0 \ \ \ \mbox{for} \ r > t,
\end{align}
i.e. $\mathbb{E}[D^1 \int_0^t B_{s}^{H,i_1} \delta B_s^{H,i_2}]=\mathbb{E}[D \int_0^t B_{s}^{H,i_1} \delta B_s^{H,i_2}]=0$. 
Next, we compute $\frac{1}{2}\mathbb{E}[D^2 \int_0^t B_{s}^{H,i_1} \delta B_s^{H,i_2}]$. The second Malliavin derivative $D^2 \int_0^t B_{s}^{H,i_1} \delta B_s^{H,i_2}$ is given as follows: for $t_1,t_2\leq t$, 
\begin{eqnarray}
&&D_{\ell_1,t_1}D_{\ell_2,t_2}\int_0^t B_{s}^{H,i_1} \delta B_s^{H,i_2}\nn\\
&=&D_{\ell_1,t_1} (B_{t_2}^{H,i_1} {\bf 1}_{\ell_2=i_2} + (B_{t}^{H,i_2}-B_{t_2}^{H,i_2}) {\bf 1}_{\ell_2=i_1} )\nn\\
&=& {\bf 1}_{0\leq t_1\leq t_2 \leq t}{\bf 1}_{\ell_1=i_1} {\bf 1}_{\ell_2=i_2}+{\bf 1}_{0\leq t_1,t_2 \leq t}{\bf 1}_{\ell_1=i_2} {\bf 1}_{\ell_2=i_1}-{\bf 1}_{0\leq t_1\leq t_2 \leq t}{\bf 1}_{\ell_1=i_2} {\bf 1}_{\ell_2=i_1}, \ \ \ell_1,\ell_2=1,\ldots,d,\nn\\
\end{eqnarray}
from (\ref{first_mall_deriv}). Then, for $\ell_1,\ell_2=1,\ldots,d$, the map $(t_1,t_2) \mapsto  \mathbb{E}[D_{\ell_1,t_1}D_{\ell_2,t_2}\int_0^t B_{s}^{H,i_1} \delta B_s^{H,i_2}]$ has the representation: 
\begin{align}
(t_1,t_2) \mapsto  & \ \ \mathbb{E}[D_{\ell_1,t_1}D_{\ell_2,t_2}\int_0^t B_{s}^{H,i_1} \delta B_s^{H,i_2}] \nn\\
&={\bf 1}_{[0,t_2]}(t_1) {\bf 1}_{[0,t]}(t_2) {\bf 1}_{\ell_1=i_1} {\bf 1}_{\ell_2=i_2} +({\bf 1}_{[0,t]}(t_1) -{\bf 1}_{[0,t_2]}(t_1)) {\bf 1}_{[0,t]}(t_2) {\bf 1}_{\ell_1=i_2} {\bf 1}_{\ell_2=i_1},
\end{align}
and for ${\bf 1}_{[0,t]}e_{\ell_1},{\bf 1}_{[0,t]}e_{\ell_2} \in \bar{\cal H}^d$, 
we have 
\begin{align}
&\Big\langle {\bf 1}_{[0,t]}e_{\ell_1} \otimes {\bf 1}_{[0,t]}e_{\ell_2}, \mathbb{E}[D^2 \int_0^t B_{s}^{H,i_1} \delta B_s^{H,i_2}] \Big\rangle_{(\bar{\cal H}^d)^{\otimes 2}} \nn \\
=& 
\int_0^t \int_0^t 
K_H(t,t_1){\bf 1}_{[0,t]}(t_1) 
K_H(t,t_2){\bf 1}_{[0,t]}(t_2) 
\Big\{ K_H(t_2,t_1){\bf 1}_{[0,t_2]}(t_1) K_H(t,t_2) {\bf 1}_{[0,t]}(t_2) {\bf 1}_{\ell_1=i_1} {\bf 1}_{\ell_2=i_2} \nn \\
& ~~~~~~~~~~~~~~~  + \{ K_H(t,t_1) {\bf 1}_{[0,t]}(t_1) - K_H(t_2,t_1) {\bf 1}_{[0,t_2]}(t_1)\} K_H(t,t_2) {\bf 1}_{[0,t]}(t_2) {\bf 1}_{\ell_1=i_2} {\bf 1}_{\ell_2=i_1} \Big\} dt_1dt_2.
\end{align}
We note
\begin{align}
&\frac{1}{p!}\mathbb{E}[D^p \int_0^t B_{s}^{H,i_1} \delta B_s^{H,i_2}]= 0, \ \ p\geq 3. 
\end{align}

Thus, we only need to compute the following term: 
\begin{align}
&\sum_{\gamma_1,\gamma_2=1}^d \mathbb{E} [ \partial_{\gamma_1}\partial_{\gamma_2}\delta_y (B_t^H) \langle DB_t^{H,\gamma_1} \otimes DB_t^{H,\gamma_2}, \frac{1}{2!} \mathbb{E}[D^2 \int_0^t B_{s}^{H,i_1} \delta B_s^{H,i_2}] \rangle_{(\bar{\cal H}^d)^{\otimes 2}} ].
\end{align}
Since $D_{\ell,\cdot} B_t^{H,i}={\bf 1}_{[0,t]}(\cdot) {\bf 1}_{\ell=i} $ for $i,\ell=1,\ldots,d$, we have
\begin{align}
&\sum_{\gamma_1,\gamma_2=1}^d \mathbb{E} [ \partial_{\gamma_1}\partial_{\gamma_2}\delta_y (B_t^H) \langle DB_t^{H,\gamma_1} \otimes DB_t^{H,\gamma_2}, \frac{1}{2!} \mathbb{E} [D^2 \int_0^t B_{s}^{H,i_1} \delta B_s^{H,i_2}] \rangle_{(\bar{\cal H}^d)^{\otimes 2}} ]  \nn\\
=&
\sum_{\gamma_1,\gamma_2=1}^d \mathbb{E} [ \partial_{\gamma_1}\partial_{\gamma_2}\delta_y (B_t^H) \frac{1}{2!} \sum_{\ell_1,\ell_2=1}^d \int_0^t \int_0^t K_H(t,t_1) {\bf 1}_{[0,t]}(t_1) {\bf 1}_{\ell_1=\gamma_1} K_H(t,t_2) {\bf 1}_{[0,t]}(t_2) {\bf 1}_{\ell_2=\gamma_2}  \nn \\
& \hspace{11.25em} K_H(t,t_2) {\bf 1}_{[0,t]}(t_2) \{ K_H(t_2,t_1){\bf 1}_{[0,t_2]}(t_1) {\bf 1}_{\ell_1=i_1} {\bf 1}_{\ell_2=i_2}  \nn\\
&\hspace{11.5em} +\{K_H(t,t_1) {\bf 1}_{[0,t]}(t_1) - K_H(t_2,t_1) {\bf 1}_{[0,t_2]}(t_1)\}{\bf 1}_{\ell_1=i_2} {\bf 1}_{\ell_2=i_1} \} dt_1 dt_2 ]  \nn\\
=&
\mathbb{E} [ \partial_{i_1}\partial_{i_2}\delta_y (B_t^H)] \frac{1}{2} \int_0^t \int_0^t K_H(t,t_1)  K_H(t,t_2)^2 K_H(t_2,t_1){\bf 1}_{[0,t_2]}(t_1) dt_1 dt_2 \nn \\
&+
\mathbb{E} [ \partial_{i_2}\partial_{i_1}\delta_y (B_t^H)] \frac{1}{2} \int_0^t \int_0^t K_H(t,t_1)  K_H(t,t_2)^2  \{ K_H(t,t_1)  - K_H(t_2,t_1) {\bf 1}_{[0,t_2]}(t_1) \} dt_1dt_2.
\end{align}
Here, we note that
\begin{align}
&\mathbb{E} [ \partial_{i_1}\partial_{i_2}\delta_y (B_t^H)] 
=
\mathbb{E} [ \partial_{i_2}\partial_{i_1}\delta_y (B_t^H)] \nn\\
=&
\mathbb{E} [ \delta_y (B_t^H) H_{(i_1,i_2)}(B_t^H,1) ] \ (=\mathbb{E} [ \delta_y (B_t^H) H_{(i_2,i_1)}(B_t^H,1) ]) \nn \\
=&
\mathbb{E} [ \delta_y (B_t^H) \frac{1}{t^{4H}} \{ B_t^{H,i_1}B_t^{H,i_2}- t^{2H} {\bf 1}_{i_1=i_2} \} ].
\end{align}
Therefore, we obtain
\begin{align}
&\sum_{\gamma_1,\gamma_2=1}^d \mathbb{E} [ \partial_{\gamma_1}\partial_{\gamma_2}\delta_y (B_t^H) \langle DB_t^{H,\gamma_1} \otimes DB_t^{H,\gamma_2}, \frac{1}{2!} \mathbb{E} [D^2 \int_0^t B_{s}^{H,i_1} \delta B_s^{H,i_2}] \rangle_{(\bar{\cal H}^d)^{\otimes 2}} ] \nn\\
=&
\mathbb{E} [ \partial_{i_1}\partial_{i_2}\delta_y (B_t^H)] \frac{1}{2} \int_0^t \int_0^t K_H(t,t_1)  K_H(t,t_2)^2  K_H(t,t_1) dt_1dt_2 \nn\\
=&
\mathbb{E} [ \delta_y (B_t^H) \frac{1}{t^{4H}}\{ B_t^{H,i_1}B_t^{H,i_2}- t^{2H} {\bf 1}_{i_1=i_2} \} ]
\frac{1}{2} \int_0^t K_H(t,t_1)^2 dt_1 \int_0^t   K_H(t,t_2)^2 dt_2 \nn\\
=&
\mathbb{E} [ \delta_y (B_t^H) \frac{1}{2} \{ B_t^{H,i_1}B_t^{H,i_2}- t^{2H} {\bf 1}_{i_1=i_2} \} ],
\end{align}
since $R_H(t,t)=\int_0^t K_H(t,s)^2 ds=t^{2H}$. Then (\ref{exp_iterated_integral}) becomes 
\begin{eqnarray}
&&\mathbb{E} [\delta_y(B_t^{H})\int_0^t \int_0^{t_{2}} \circ dB_{t_1}^{H,i_1} \circ dB_{t_2}^{H,i_2}]\nn\\
&=&
\mathbb{E} [ \delta_y(B_t^{H}) \frac{1}{2} \{ B_t^{H,i_1}B_t^{H,i_2} - t^{2H} {\bf 1}_{i_1=i_2} \} ]  +\mathbb{E}[\delta_y(B_t^{H}) ] \frac{1}{2} t^{2H}{\bf 1}_{i_1=i_2} \nn\\
&=&
\mathbb{E} [ \delta_y(B_t^{H}) \frac{1}{2} B_t^{H,i_1}B_t^{H,i_2}  ], 
\end{eqnarray}
and for $f \in {\cal S}'(\mathbb{R}^e)$, we have
\begin{eqnarray}
&&t^{2H}\mathbb{E}[ \partial_{j} f ( \bar{X}_1^{x,t^H} ) V_{i_1}V^{j}_{i_2}(x) \mathbb{B}^{H}_{(i_1,i_2)}(1) ] \nn\\
&=&t^{2H}\mathbb{E}[ \partial_{j} f ( \bar{X}_1^{x,t^H} ) \frac{1}{2}V_{i_1}V^{j}_{i_2}(x) B_1^{H,i_1}B_1^{H,i_2}  ] \ \ \mbox{for} \ \ j=1,\ldots,e. 
\end{eqnarray}
Finally, we have the following by the integration by parts: 
\begin{align}
&t^{2H}\mathbb{E}[ \partial_{\alpha_1} f ( \bar{X}_1^{x,t^H} ) \frac{1}{2} V_{i_1}V^{\alpha_1}_{i_2}(x) B_1^{H,i_1}B_1^{H,i_2}  ]
=t^{2H} \mathbb{E}[f ( \bar{X}_1^{x,t^H} ) H_{(\alpha_1)} (\bar{X}_1^{x,t^H}, \frac{1}{2}V_{i_1}V^{\alpha_1}_{i_2}(x) B_1^{H,i_1}B_1^{H,i_2})] \nn\\
=&
t^{2H}\mathbb{E} [ f ( \bar{X}_1^{x,t^H} ) \sum_{\alpha_2=1}^e \sum_{i_3=1}^d V_{i_1}V^{\alpha_1}_{i_2}(x) \nn\\
& \hspace{5em} t^H V_{i_3}^{\alpha_2}(x) A^{-1}_{\alpha_1,\alpha_2}(x) \frac{1}{2t^{2H}}\{ B_1^{H,i_1}B_1^{H,i_2}B_1^{H,i_3} - B_1^{H,i_1} {\bf 1}_{i_2=i_3\neq 0} - B_1^{H,i_2} {\bf 1}_{i_1=i_3\neq 0}  \} ] \nn\\
=&
t^{H}\mathbb{E} [ f ( \bar{X}_1^{x,t^H} ) \sum_{\alpha_2=1}^e \sum_{i_3=1}^d V_{i_1}V^{\alpha_1}_{i_2}(x) \nn\\
& \hspace{5em} V_{i_3}^{\alpha_2}(x) A^{-1}_{\alpha_1,\alpha_2}(x) \frac{1}{2}\{ B_1^{H,i_1}B_1^{H,i_2}B_1^{H,i_3} - B_1^{H,i_1} {\bf 1}_{i_2=i_3\neq 0} - B_1^{H,i_2} {\bf 1}_{i_1=i_3\neq 0}  \} ].
\end{align}
Since the formula holds for all $f \in {\cal S}'(\mathbb{R}^e)$, we have the expansion of the probability distribution function with the error order $O(t^{1-H})$. Moreover, since we have the following integration by parts: for $f \in {\cal S}(\mathbb{R}^e)$ and $\alpha_1=1,\ldots,e$, 
\begin{align}
t \mathbb{E}[ \partial_{\alpha_1} f ( \bar{X}_1^{x,t^H} )  V^{\alpha_1}_{0}(x)  ]
=&
t \mathbb{E} [ f ( \bar{X}_1^{x,t^H} ) \sum_{\alpha_2=1}^e \sum_{i_1=1}^d V^{\alpha_1}_{0}(x) t^H V_{i_1}^{\alpha_2}(x) A^{-1}_{\alpha_1,\alpha_2}(x) \frac{1}{t^{2H}}B_1^{H,i_1}  ]\nn\\
=&
t^{1-H} \mathbb{E} [ f ( \bar{X}_1^{x,t^H} ) \sum_{\alpha_2=1}^e \sum_{i_1=1}^d V^{\alpha_1}_{0}(x) V_{i_1}^{\alpha_2}(x) A^{-1}_{\alpha_1,\alpha_2}(x) B_1^{H,i_1}  ],
\end{align}
the desired $O(t^{2H})$-expansion is obtained. $\Box$ \\ \\

We finally show a numerical example for the expansion of Theorem \ref{expansion_rde} in order to validate the result. Consider the following rough differential equation driven by fractional Brownian motion: 
\begin{align}
dX_t^x=V(X_t^x)\circ d{\bf B}_t^H, \ \ X_0^x=x \in \mathbb{R},
\end{align}
where $V:x \mapsto V(x)=\sigma x$ with $\sigma>0$. We compute the probability distribution function $z \mapsto \mathbb{P}(X_t^x \leq z)$ using the $O(t^{2H})$-asymptotic expansion of Theorem \ref{expansion_rde}. We set the parameters as $H=0.4$, $t=0.25$, $\sigma=0.3$, $x=10$, $z\in [5,15]$. We compare the asymptotic expansion with the normal approximation $\mathbb{P}(X_t^x \leq z) \approx \mathbb{P}(\bar{X}_1^{x,t^H} \leq z)$, and the exact solution. The asymptotic expansion and the normal approximation are implemented by quasi-Monte Carlo method with $10^6$-paths. 

The following figure (Figure 1) shows the effectiveness of our asymptotic expansion as its accuracy overcomes that of 
normal approximation. 

\begin{figure}[H]
\begin{center}
\includegraphics[scale=0.55]{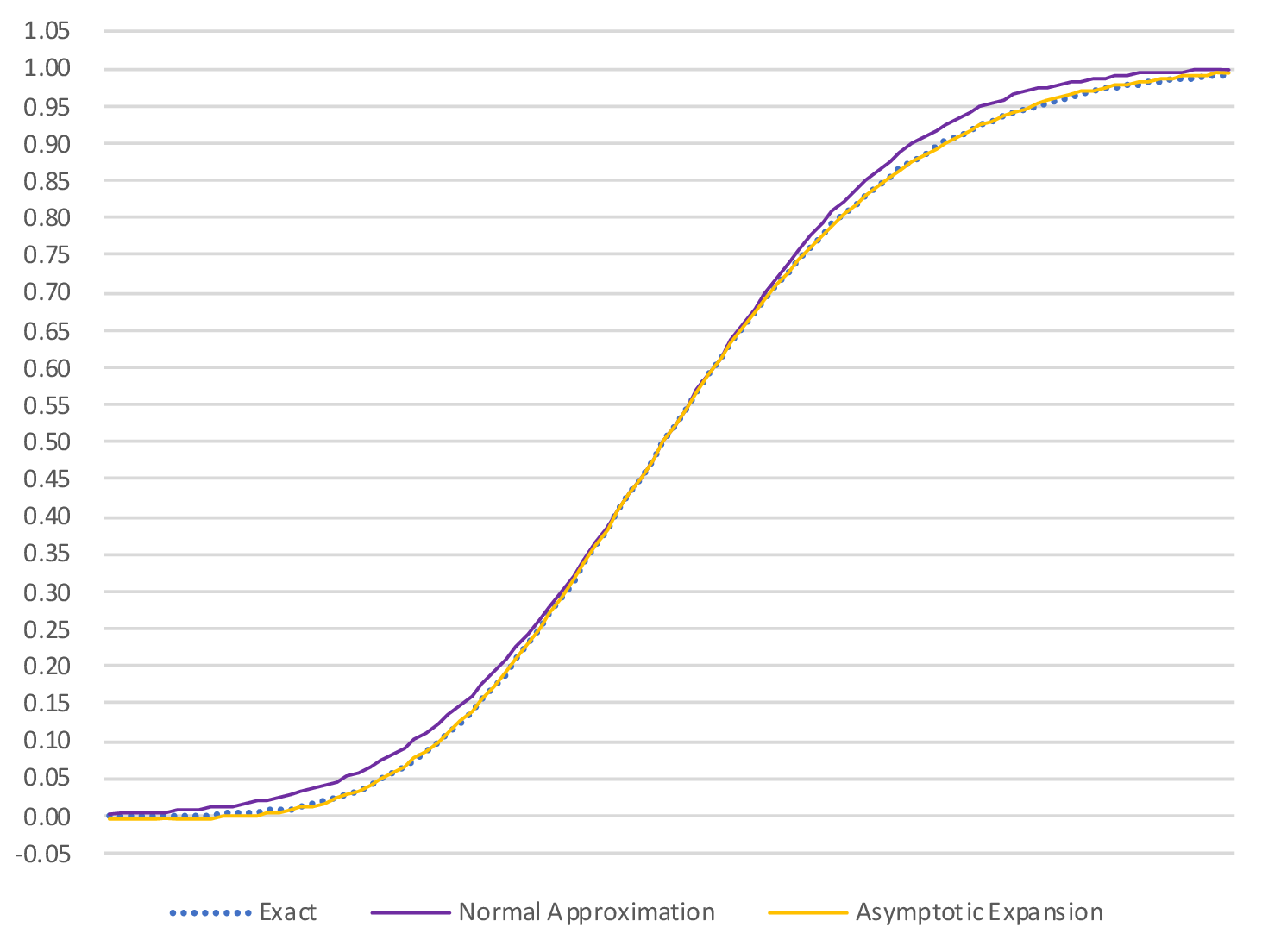}
\end{center}
\caption{Accuracy of asymptotic expansion for probability distribution function of solution to rough differential equation driven by fractional differential equation with Hurst index $H=0.4$}
\end{figure}

\section{Concluding remarks}
In the paper, we have provided a new asymptotic expansion formula of expectation of general multidimensional Wiener functionals. The uniform estimate of the asymptotic expansion has been obtained under a weaker condition on the Malliavin covariance matrix of the target Wiener functionals. Then we have shown a tractable expansion for the expectation of an irregular functional of the solution to a multidimensional RDE driven by fractional Brownian motion with Hurst index $H<1/2$. The result has been justified by a numerical example for the asymptotic expansion of a probability distribution function through a comparison with the normal approximation. We note that it is possible to give the similar expansion for expectations of multidimensional differential equations driven by fractional Brownian motion with Hurst index $H\in (1/2,1)$.

It is interesting to see whether the proposed asymptotic expansion approach can be applied to numerical methods such as Monte-Carlo methods or discretization methods for RDEs as in Takahashi and Yoshida (2005) \cite{TYo}, Takahashi and Yamada (2016) \cite{TY16} and Yamada (2019) \cite{Y19} for standard SDEs, which will be studied as future works. 

\section*{Acknowledgements}
We are grateful to the associate editor and the two anonymous reviewers for their valuable comments and suggestions. This work is supported by JST PRESTO (Grant Number JPMJPR2029), Japan.

\end{document}